\documentclass{icm2010}
\usepackage{pstricks}
\usepackage{pst-node}
\usepackage{amssymb}
\newtheorem{theorem}{Theorem}
\newtheorem{lemma}[theorem]{Lemma}
\newtheorem{proposition}[theorem]{Proposition}

\newtheorem{conjecture}[theorem]{Conjecture}
\newtheorem{defi}[theorem]{Definition}
\theoremstyle{remark}

\newtheorem*{remark*}{Remark}

\newtheorem*{example*}{Example}
\newcommand{\Z}{\mathbb{Z}}
\newcommand{\R}{\mathbb{R}}
\newcommand{\C}{\mathbb{C}}

\newcommand{\Sym}{\mathrm{Sym}}
\renewcommand{\Re}{\mathrm{Re}\,}
\renewcommand{\Im}{\mathrm{Im}\,}
\newcommand{\F}{\mathcal{F}}
\newcommand{\FF}{\mathbb{F}}
\newcommand{\A}{\mathcal{A}}
\newcommand{\Y}{\mathcal{Y}}
\newcommand{\T}{\mathbb{T}}
\renewcommand{\SS}{\mathcal{S}}
\author{Denis Auroux}
\contact[auroux@math.berkeley.edu]{Department of Mathematics, UC Berkeley,
Berkeley CA 94720-3840, USA}
\title{Fukaya categories and bordered Heegaard-Floer homology}
\begin{document}
\begin{abstract}
We outline an interpretation of Heegaard-Floer homology of 3-manifolds
(closed or with boundary)
in terms of the symplectic topology
of symmetric products of Riemann surfaces, as suggested by recent
work of Tim Perutz and Yank\i{} Lekili. In particular we discuss the
connection between the Fukaya category of the symmetric product and the
bordered algebra introduced by Robert Lipshitz, Peter Ozsv\'ath and Dylan Thurston, and
recast bordered Heegaard-Floer homology in this language.
\end{abstract}



\maketitle


\section{Introduction}

In its simplest incarnation, Heegaard-Floer homology associates to a
closed 3-manifold $Y$ a graded abelian group $\widehat{HF}(Y)$. This
invariant is constructed
by considering a Heegaard splitting $Y=Y_1\cup_\Sigma Y_2$ of $Y$ into two
genus $g$ handlebodies, each of which determines a product torus in the
$g$-fold symmetric product of the Heegaard surface $\Sigma=\partial Y_1=
-\partial Y_2$. Deleting a marked point $z$ from $\Sigma$ to obtain an open
surface, $\widehat{HF}(Y)$ is then defined as the 
Lagrangian Floer homology of the two tori $T_1,T_2$
in $\Sym^g(\Sigma\setminus \{z\})$, see \cite{OS}.

It is natural to ask how more general decompositions of 3-manifolds 
fit into this picture, and whether Heegaard-Floer theory can be viewed
as a TQFT (at least in some partial sense). From the point of view of
symplectic geometry, a natural answer is suggested
by the work of Tim Perutz and Yank\i{} Lekili. Namely,
an elementary cobordism between two connected Riemann surfaces 
$\Sigma_1,\Sigma_2$ given by attaching a single handle determines a 
Lagrangian correspondence between appropriate symmetric products of
$\Sigma_1$ and $\Sigma_2$~\cite{Perutz}. By composing these correspondences,
one can associate to a 3-manifold with connected boundary $\Sigma$
of genus $g$
a {\em generalized Lagrangian submanifold} (cf.\ \cite{WW}) of the symmetric product
$\Sym^g(\Sigma\setminus \{z\})$.
Recent work of Lekili and Perutz \cite{LP} shows that, given a decomposition 
$Y=Y_1\cup_\Sigma Y_2$ of a closed 3-manifold, the (quilted) 
Floer homology of the generalized Lagrangian submanifolds of
$\Sym^g(\Sigma\setminus \{z\})$ determined by $Y_1$ and $Y_2$ recovers
$\widehat{HF}(Y)$.

From a more combinatorial perspective, the {\it bordered Heegaard-Floer
homology} of 
Robert Lipshitz, Peter Ozsv\'ath and Dylan Thurston \cite{LOT} associates
to a parameterized Riemann surface $F$ with connected boundary a finite
dimensional 
differential algebra $\A(F)$ over $\Z_2$, and to a 3-manifold $Y$ with boundary
$\partial Y= F\cup D^2$ a right $A_\infty$-module $\widehat{CFA}(Y)$
over $\A(F)$, as well as a left dg-module $\widehat{CFD}(Y)$ over $\A(-F)$.
The main result of \cite{LOT} shows that, given a decomposition of a closed 
3-manifold $Y=Y_1\cup Y_2$ with $\partial Y_1=-\partial Y_2= F\cup D^2$, 
$\widehat{HF}(Y)$ can be computed in terms of the $A_\infty$-tensor product
of the modules associated to $Y_1$ and $Y_2$, namely
$$\widehat{HF}(Y)\simeq H_*(\widehat{CFA}(Y_1)\otimes_{\A(F)}
\widehat{CFD}(Y_2)).$$

In order to connect these two approaches, we consider a {\it partially
wrapped} version of Floer theory for product Lagrangians in
symmetric products of open Riemann surfaces. Concretely, given a Riemann
surface with boundary $F$, a finite collection $Z$ of marked points 
on $\partial F$, and an integer $k\ge 0$, we consider a partially 
wrapped Fukaya category $\F(\Sym^k(F),Z)$, which differs from the usual
(compactly supported) Fukaya category by the inclusion of additional
objects, namely products of disjoint properly embedded
arcs in $F$ with boundary in $\partial F\setminus Z$. A nice feature of
these categories is that they admit explicit sets of 
generating objects:

\begin{theorem} \label{thm:generate}
Let $F$ be a compact Riemann surface with non-empty boundary, $Z$ a finite
subset of $\partial F$, and
$\underline{\alpha}=\{\alpha_1,\dots,\alpha_n\}$ a collection of
disjoint properly embedded arcs in $F$ with boundary in $\partial F\setminus
Z$. Assume that $F\setminus (\alpha_1\cup\dots\cup \alpha_n)$ is a union of
discs, each of which contains at most one point of $Z$. Then for $0\le k\le n$,
the partially wrapped Fukaya category $\F(\Sym^k(F),Z)$ is generated by the
$\binom{n}{k}$ Lagrangian submanifolds $D_s=\prod_{i\in s}\alpha_i$,
where $s$ ranges over all $k$-element subsets of $\{1,\dots,n\}$.
\end{theorem}

To a decorated surface $\FF=(F,Z,\underline{\alpha}=\{\alpha_i\})$ we can
associate an $A_\infty$-algebra 
\begin{equation}\label{eq:AFk}\textstyle
\A(\FF,k)=\bigoplus\limits_{s,t} \hom_{\F(\Sym^k(F),Z)}(D_s,D_t).
\end{equation}
The following special case is of particular interest:

\begin{theorem} \label{thm:algebra}
Assume that $F$ has a single boundary component, $|Z|=1$,
and the arcs $\alpha_1,\dots,\alpha_n$ $(n=2g(F))$ decompose $F$
into a single disc. Then
$\A(\FF,k)$ coincides with Lipshitz-Ozsv\'ath-Thurston's {\em bordered
algebra \cite{LOT}}.
\end{theorem}

\noindent
(The result remains true in greater generality, the only key
requirement being that every component of $F\setminus (\alpha_1\cup
\dots\cup \alpha_n)$ should contain at least one point of $Z$.)
\medskip

Now, consider a {\em sutured} 3-manifold $Y$, i.e.\ a 3-manifold $Y$
with non-empty boundary, equipped with a decomposition $\partial Y=
(-F_-)\cup F_+$, where $F_\pm$ are oriented surfaces with boundary.
Assume moreover that $\partial Y$ and $F_\pm$ are connected, and
denote by $g_\pm$ the genus of $F_\pm$. Given two integers
$k_\pm$ such that $k_+-k_-=g_+-g_-$ and
a suitable Morse function
on $Y$, Perutz's construction associates to $Y$ a generalized Lagrangian
correspondence (i.e.\ a formal composition of Lagrangian correspondences)
$\mathbb{T}_Y$ from $\Sym^{k_-}(F_-)$ to $\Sym^{k_+}(F_+)$.
By the main result of \cite{LP} this
correspondence is essentially independent of the chosen
Morse function. 

Picking a finite set of marked points $Z\subset \partial F_-=\partial F_+$,
and two collections
of disjoint arcs $\underline{\alpha}_-$ and $\underline{\alpha}_+$ on $F_-$ and 
$F_+$, we have two decorated surfaces
$\FF_\pm=(F_\pm,Z,\underline{\alpha}_\pm)$, and collections of product
Lagrangian submanifolds $D_{\pm,s}$ \hbox{($s\in \mathcal{S}_\pm$)} in 
$\Sym^{k_\pm}(F_\pm)$  (namely, all products of $k_\pm$ of the arcs 
in $\underline{\alpha}_\pm$). By a Yoneda-style construction, the correspondence $\mathbb{T}_Y$ then
determines an $A_\infty$-bimodule
\begin{equation}\label{eq:YTY}
\textstyle \mathcal{Y}(\mathbb{T}_Y)=\bigoplus\limits_{(s,t)\in \mathcal{S}_-
\times\mathcal{S}_+} \hom(D_{-,s},\mathbb{T}_Y,D_{+,t})\in
\A(\FF_-,k_-)\text{-mod-}\A(\FF_+,k_+),\end{equation} where
$\hom(D_{-,s},\mathbb{T}_Y,D_{+,t})$ is defined 
in terms of quilted Floer complexes \cite{MWW,WW,WW2} after suitably
perturbing $D_{-,s}$ and $D_{+,t}$ by partial wrapping along the boundary. A slightly
different but equivalent definition is as follows.
With quite a bit of extra work, via the Ma'u-Wehrheim-Woodward machinery
the correspondence $\mathbb{T}_Y$ defines an
$A_\infty$-functor $\Phi_Y$ from $\F(\Sym^{k_-}(F_-),Z)$ to a suitable
enlargement of $\F(\Sym^{k_+}(F_+),Z)$; with this understood,
$\mathcal{Y}(\mathbb{T}_Y)\simeq \bigoplus_{(s,t)} \hom(\Phi_Y(D_{-,s}),D_{+,t})$.

The $A_\infty$-bimodules $\mathcal{Y}(
\mathbb{T}_Y)$ are expected to obey the following
gluing property:

\begin{conjecture} \label{conj:pairing}
Let $F,F',F''$ be
connected Riemann surfaces and $Z$ a finite subset of 
$\partial F\simeq \partial F'\simeq
\partial F''$. Let $Y_1,Y_2$ be two sutured manifolds with
$\partial Y_1=(-F)\cup F'$ and $\partial Y_2=(-F')\cup F''$, and let
$Y=Y_1\cup_{F'} Y_2$ be the sutured manifold obtained by gluing $Y_1$ and
$Y_2$ along $F'$. Equip $F,F',F''$ with collections of disjoint properly
embedded arcs
$\underline{\alpha},\underline{\alpha}',\underline{\alpha}''$, and assume
that $\underline{\alpha}'$ decomposes $F'$ into a union of discs each
containing at most one point of $Z$. Then
\begin{equation}\label{eq:pairing}
\mathcal{Y}(\mathbb{T}_{Y})\simeq \mathcal{Y}(\mathbb{T}_{Y_1})
\otimes_{\A(\FF',k')} \mathcal{Y}(\mathbb{T}_{Y_2}).
\end{equation}
\end{conjecture}

In its most general form this statement relies on
results in Floer theory for generalized Lagrangian
correspondences which are not yet fully established, hence we
state it as a conjecture; however, we believe that a proof should
be within reach of standard techniques.

As a special case, let $F$ be a genus $g$ surface with connected boundary, decorated
with a single point $z\in \partial F$ and a collection of $2g$ arcs cutting
$F$ into a disc. Then to a 3-manifold $Y_1$ with boundary $\partial Y_1=F\cup
D^2$ we can associate a generalized Lagrangian submanifold $\mathbb{T}_{Y_1}$
of $\Sym^g(F)$, and an $A_\infty$-module $\mathcal{Y}(\mathbb{T}_{Y_1})=
\bigoplus_s \hom(\mathbb{T}_{Y_1},D_s)\in \text{mod-}\A(\FF,g)$. Viewing
$\mathbb{T}_{Y_1}$ as a generalized correspondence from $\Sym^g(-F)$ to 
$\Sym^0(D^2)=\{\text{pt}\}$
instead, we obtain a left $A_\infty$-module over $\A(-\FF,g)$. However,
$\A(-\FF,g)=\A(\FF,g)^{op}$, and the two constructions yield the same module.
If now we have another 3-manifold $Y_2$ with $\partial Y_2=-F\cup D^2$, we
can associate to it a generalized Lagrangian submanifold $\mathbb{T}_{Y_2}$
in $\Sym^g(-F)$ or, after orientation reversal, $\mathbb{T}_{-Y_2}$ in
$\Sym^g(F)$. This yields $A_\infty$-modules $\mathcal{Y}(\mathbb{T}_{Y_2})
\in \text{mod-}\A(-\FF,g)\simeq \A(\FF,g)\text{-mod}$, and
$\mathcal{Y}(\mathbb{T}_{-Y_2}) \in \text{mod-}\A(\FF,g)$. 

\begin{theorem}\label{cor:pairing2}
With this
understood, and denoting by $Y$ the closed 3-manifold obtained by gluing
$Y_1$ and $Y_2$ along their boundaries, we have quasi-isomorphisms
\begin{align}\label{eq:pairing2}\nonumber
\widehat{CF}(Y)\simeq \hom_{\F^\#(\Sym^g(F))}(\T_{Y_1},
\T_{-Y_2})&\simeq \hom_{\text{mod-}\A(\FF,g)}(\Y(\T_{-Y_2}),\Y(\T_{Y_1}))\\
&\simeq \Y(\T_{Y_1})\otimes_{\A(\FF,g)} \Y(\T_{Y_2}).
\end{align}
\end{theorem}

\noindent
In fact, $\Y(\T_{Y_i})$ is quasi-isomorphic to the bordered
$A_\infty$-module $\widehat{CFA}(Y_i)$. In light of this, it is
instructive to compare Theorem \ref{cor:pairing2} with the pairing
theorem obtained by Lipshitz, Ozsv\'ath and Thurston \cite{LOT}: 
even though $\widehat{CFA}(Y_i)$ and $\widehat{CFD}(Y_i)$
seem very different at first glance (and even at second glance), our
result suggests that they can in fact be used interchangeably.
\medskip

The rest of this paper is structured as follows: first, in section \ref{s:LP} we
explain how Heegaard-Floer homology can be understood in terms of
Lagrangian correspondences, following the work of Perutz and Lekili
\cite{Perutz,LP}. Then in section \ref{s:fukaya} we introduce partially
wrapped Fukaya categories of symmetric products, and sketch the proofs
of Theorems \ref{thm:generate} and \ref{thm:algebra}. In section
\ref{s:yoneda} we briefly discuss Yoneda embedding as well as Conjecture
\ref{conj:pairing} and Theorem \ref{cor:pairing2}. Finally, in section
\ref{s:bordered} we discuss the relation with bordered Heegaard-Floer
homology.

The reader will not find detailed proofs for any of the statements here,
nor a general discussion of partially wrapped Fukaya categories. Some of
the material is treated in greater depth in the preprint \cite{fuksymg}, 
the rest will appear in a future paper.

%

\subsection*{Acknowledgements}
I am very grateful to Mohammed Abouzaid, Sheel Ganatra, Yank\i{} Lekili,
Robert Lipshitz, Peter Ozsv\'ath,
Tim Perutz, Paul Seidel and Dylan Thurston for many helpful discussions.
I would also like to thank Ivan Smith for useful comments on the exposition.
This work was partially supported by NSF grants DMS-0600148 and DMS-0652630.

\section{Heegaard-Floer homology from Lagrangian correspondences}
\label{s:LP}

\subsection{Lagrangian correspondences}
A {\it Lagrangian correspondence} between two symplectic manifolds 
$(M_1,\omega_1)$ and $(M_2,\omega_2)$ is, by definition, a Lagrangian submanifold
of the product $M_1\times M_2$ equipped with the product symplectic 
form $(-\omega_1)\oplus \omega_2$. Lagrangian correspondences can be 
thought of as a far-reaching generalization of symplectomorphisms 
(whose graphs are examples of correspondences); in particular, under
suitable transversality assumptions we can consider the {\it composition}
of two correspondences $L_{01}\subset M_0\times M_1$ and $L_{12}\subset
M_1\times M_2$, $$L_{01}\circ L_{12}=\{(x,z)\in M_0\times M_2\,|\,
\exists\, y\in M_1\ \text{s.t.}\ (x,y)\in L_{01}\ \text{and}\ (y,z)\in
L_{12}\}.$$
The image of a Lagrangian submanifold $L_1\subset M_1$ by
a Lagrangian correspondence $L_{12}\subset M_1\times M_2$ is defined
similarly, viewing $L_1$ as a correspondence from $\{pt\}$ to~$M_1$.
Unfortunately, in general the geometric composition is not a smooth embedded
Lagrangian. Nonetheless, we can enlarge symplectic geometry by
considering {\it generalized Lagrangian correspondences}, i.e.\
sequences of Lagrangian correspondences (interpreted as formal
compositions), and {\it generalized Lagrangian submanifolds}, i.e.\
generalized Lagrangian correspondences from $\{pt\}$ to a given symplectic
manifold.

The work of Ma'u, Wehrheim and Woodward (see e.g.\
\cite{WW,WW2,MWW}) shows that Lagrangian
Floer theory behaves well with respect to (generalized) correspondences. 
Given a sequence of
Lagrangian correspondences $L_{i-1,i}\subset M_{i-1}\times M_i$
($i=1,\dots,n$), with $M_0=M_n=\{pt\}$, the {\it quilted Floer
complex} $CF(L_{0,1},\dots,L_{n-1,n})$ is generated by {\it generalized
intersections}, i.e.\ tuples $(x_1,\dots,x_{n-1})\in M_1\times \dots\times
M_{n-1}$ such that $(x_{i-1},x_i)\in L_{i-1,i}$ for all $i$, and carries
a differential which counts ``quilted pseudoholomorphic strips'' in
$M_1\times\dots\times M_{n-1}$. Under suitable technical assumptions
(e.g., monotonicity), Lagrangian Floer theory carries over to this setting.

Thus, Ma'u, Wehrheim and Woodward associate to a monotone symplectic manifold
$(M,\omega)$ its {\it extended Fukaya category} $\F^\#(M)$, whose objects
are monotone generalized Lagrangian submanifolds and in which morphisms 
are given by quilted Floer complexes. Composition of morphisms is defined by counting
quilted pseudoholomorphic discs, and as in usual Floer theory, it is only
associative up to homotopy, so $\F^\#(M)$ is an $A_\infty$-category. The
key property of these extended Fukaya categories is that a monotone
(generalized) Lagrangian correspondence $L_{12}$ from $M_1$ to $M_2$ induces
an $A_\infty$-functor from $\F^\#(M_1)$ to $\F^\#(M_2)$, which on the level
of objects is simply concatenation with $L_{12}$. Moreover, composition of
Lagrangian correspondences matches with composition of $A_\infty$-functors
\cite{MWW}.

\begin{remark*}
By construction, the usual Fukaya category $\F(M)$ admits a fully faithful
embedding as a subcategory of $\F^\#(M)$. In fact, $\F^{\#}(M)$ embeds 
into the category of $A_\infty$-modules over the usual Fukaya category, so
although generalized Lagrangian correspondences play an important conceptual
role in our discussion, they only enlarge the Fukaya category in a fairly
mild manner.
\end{remark*}

\subsection{Symmetric products}

As mentioned in the introduction, work in progress of Lekili and Perutz
\cite{LP}  shows that Heegaard-Floer homology can be understood
in terms of quilted Floer homology for Lagrangian 
correspondences between symmetric products. The relevant correspondences
were introduced by Perutz in his thesis \cite{Perutz}.

Let $\Sigma$ be an open Riemann surface (with infinite cylindrical ends, 
i.e., the complement of a finite set in a compact Riemann surface), 
equipped with an area form $\sigma$. We consider the symmetric product 
$\Sym^k(\Sigma)$, equipped with the product complex structure $J$, and 
a K\"ahler form $\omega$ which coincides with the product K\"ahler form
on $\Sigma^k$ away from the diagonal strata. Following
Perutz we choose $\omega$ so that its cohomology class is negatively
proportional to $c_1(T\Sym^k(\Sigma))$.

Let $\gamma$ be a non-separating simple closed 
curve on $\Sigma$, and $\Sigma_\gamma$ the surface obtained from $\Sigma$
by deleting a tubular neighborhood of $\gamma$ and gluing in two discs.
Equip $\Sigma_\gamma$ with a complex structure which agrees with that of
$\Sigma$ away from $\gamma$, and equip $\Sym^k(\Sigma)$ and
$\Sym^{k-1}(\Sigma_\gamma)$ with K\"ahler forms $\omega$ and $\omega_\gamma$
chosen as above.

\begin{theorem}[Perutz \cite{Perutz}]\label{th:perutz}
The simple closed curve $\gamma$ determines a Lagrangian correspondence
$T_\gamma$ in the product $(\Sym^{k-1}(\Sigma_\gamma)\times \Sym^k(\Sigma),
-\omega_\gamma\oplus \omega)$, canonically up to Hamiltonian isotopy.
\end{theorem}

Given $r$ disjoint simple closed curves $\gamma_1,\dots,\gamma_r$, linearly
independent in $H_1(\Sigma)$, we can consider the sequence of
correspondences that arise from successive surgeries along
$\gamma_1,\dots,\gamma_r$. The main properties of these correspondences
(see Theorem A in \cite{Perutz}) imply immediately that their composition
defines an embedded Lagrangian correspondence $T_{\gamma_1,\dots,\gamma_r}$ 
in $\Sym^{k-r}(\Sigma_{\gamma_1,\dots,\gamma_r})\times \Sym^k(\Sigma)$.

When $r=k=g(\Sigma)$, this construction yields
a Lagrangian torus in $\Sym^k(\Sigma)$, which by \cite[Lemma 3.20]{Perutz}
is smoothly isotopic to the product torus
$\gamma_1\times\dots\times\gamma_k$; 
Lekili and Perutz show that these two tori are in fact Hamiltonian isotopic
\cite{LP}.

\begin{remark*}
We are not quite in the setting considered by Ma'u, Wehrheim and
Woodward, but Floer theory remains well behaved thanks to two key properties
of the Lagrangian submanifolds under consideration: their relative $\pi_2$
is trivial (which prevents bubbling), and they are {\em balanced}.
(A Lagrangian submanifold in a monotone symplectic manifold is said to be
balanced if the holonomy of a fixed connection 1-form with curvature equal
to the symplectic form vanishes on it; this is a natural analogue of the
notion of exact Lagrangian submanifold in an exact symplectic manifold).
The balancing condition is closely related to admissibility of
Heegaard diagrams, and ensures that the symplectic area of a pseudo-holomorphic
strip connecting two given intersection points 
is determined {\it a priori} by its Maslov index (cf.\ \cite[Lemma
4.1.5]{WW2}). This property is what allows us to work over $\Z_2$ rather than
over a Novikov field. 
\end{remark*}

\subsection{Heegaard-Floer homology}\label{ss:LP}

Consider a closed 3-manifold $Y$, and a Morse function $f:Y\to\R$
(with only one minimum and one maximum, and with distinct critical 
values). Then the complement $Y'$ of a ball in $Y$ (obtained by deleting
a neighborhood of a Morse trajectory from the maximum to the minimum)
can be decomposed into a succession of elementary cobordisms $Y'_i$
($i=1,\dots,r$) between connected Riemann surfaces with boundary 
$\Sigma_0,\Sigma_1,\dots,\Sigma_r$ (where $\Sigma_0=\Sigma_r=D^2$, and
the genus increases or decreases by 1 at each step). By Theorem
\ref{th:perutz}, each $Y'_i$ determines a Lagrangian correspondence
$L_i\subset \Sym^{g_{i-1}}(\Sigma_{i-1})\times \Sym^{g_i}(\Sigma_i)$ between
the relevant symmetric products (here $g_i$ is the genus of $\Sigma_i$, and
we implicitly complete $\Sigma_i$ by attaching to it an 
infinite cylindrical end). By the work of Lekili and Perutz \cite{LP},
the quilted Floer homology of the sequence $(L_1,\dots,L_r)$ is independent
of the choice of the Morse function $f$ and isomorphic to $\widehat{HF}(Y)$. 

More generally, consider a sutured 3-manifold $Y$, i.e.\ a 3-manifold whose
boundary is decomposed into a union
$(-F_-)\cup_\Gamma  F_+$, where $F_\pm$ are connected oriented surfaces 
of genus $g_\pm$ with boundary $\partial F_-\simeq \partial F_+\simeq \Gamma$.
Shrinking $F_\pm$ slightly within $\partial Y$, it is advantageous to
think of the boundary of $Y$ as consisting actually of {\it three} pieces,
$\partial Y=(-F_-)\cup (\Gamma\times [0,1])\cup F_+$. By considering
a Morse function $f:Y\to [0,1]$ with index 1 and 2 critical points only,
with $f^{-1}(1)=F_-$ and $f^{-1}(0)=F_+$, we can view $Y$ as a succession
of elementary cobordisms between connected Riemann surfaces with boundary,
starting with $F_-$ and ending with $F_+$. As above, Perutz's construction
associates a Lagrangian correspondence to each of these elementary
cobordisms. Thus we can associate to $Y$ a generalized Lagrangian 
correspondence $\mathbb{T}_Y=\mathbb{T}_{Y,k_\pm}$ from $\Sym^{k_-}(F_-)$ to $\Sym^{k_+}(F_+)$ whenever
$k_+-k_-=g_+-g_-$. The generalized correspondence $\mathbb{T}_Y$ can
be viewed either as an object of the extended Fukaya category
$\F^\#(\Sym^{k_-}(-F_-)\times \Sym^{k_+}(F_+))$, or
as an $A_\infty$-functor from
$\F^\#(\Sym^{k_-}(F_-))$ to $\F^\#(\Sym^{k_+}(F_+))$.

\begin{theorem}[Lekili-Perutz \cite{LP}]\label{th:LP}
Up to quasi-isomorphism the object $\mathbb{T}_Y$ 
is independent of the choice of Morse function on $Y$.
\end{theorem}

Given two sutured manifolds $Y_1$ and $Y_2$ ($\partial
Y_i=(-F_{i,-})\cup F_{i,+}$) and a diffeomorphism
$\phi:F_{1,+}\to F_{2,-}$, gluing $Y_1$ and $Y_2$ by identifying the positive boundary of
$Y_1$ with the negative boundary of $Y_2$ via $\phi$
yields a new sutured manifold $Y'$. As a cobordism from $F_{1,-}$ to
$F_{2,+}$, $Y'$ is simply the concatenation of the cobordisms 
$Y_1$ and $Y_2$. Hence, the generalized Lagrangian correspondence
$\T_{Y'}$ associated to $Y'$ is just the (formal) composition of $\mathbb{T}_{Y_1}$ and
$\mathbb{T}_{Y_2}$.

The case where $Y_1$ is a cobordism from the disc $D^2$ to a genus $g$
surface $F$ (with a single boundary component) and $Y_2$ is a cobordism
from $F$ to $D^2$ (so $\partial Y_1\simeq -\partial Y_2\simeq F\cup_{S^1}
D^2$) is of particular interest. In that case, we associate to $Y_1$ a
generalized correspondence from $\Sym^0(D^2)=\{pt\}$ to $\Sym^g(F)$, i.e.\
an object $\T_{Y_1}$ of $\F^\#(\Sym^g(F))$, and to $Y_2$ a generalized
correspondence $\T_{Y_2}$ from $\Sym^g(F)$ to $\Sym^0(D^2)=\{pt\}$, i.e.\ a generalized
Lagrangian submanifold of $\Sym^g(-F)$.
Reversing the orientation of $Y_2$, i.e.\ viewing $-Y_2$ as the opposite 
cobordism from $D^2$ to $F$, we get a generalized Lagrangian submanifold
$\mathbb{T}_{-Y_2}$ in $\Sym^g(F)$, which differs from $\mathbb{T}_{Y_2}$
simply by orientation reversal. Denoting by $Y$ $(=Y'\cup B^3)$ the closed 3-manifold
obtained by gluing $Y_1$ and $Y_2$ along their entire boundary,
the result of \cite{LP} now says that 
\begin{equation}\label{eq:LPpairing}
\widehat{CF}(Y)\simeq
CF(\T_{Y_1},\T_{Y_2})\simeq \hom_{\F^\#(\Sym^g(F))}(\T_{Y_1},\T_{-Y_2}).
\end{equation}

\section{Partially wrapped Fukaya categories of symmetric products}
\label{s:fukaya}

\subsection{Positive perturbations and partial wrapping}\label{ss:wrap}
Let $F$ be a connected Riemann surface with non-empty boundary, and $Z$
a finite subset of $\partial F$. Assume for now that every connected
component of $\partial F$ contains at least one point of $Z$. Then
the components of $\partial F\setminus Z$ are open intervals, and carry a
natural orientation induced by that of $F$. 

\begin{defi}
Let $\underline\lambda=(\lambda_1,\dots,\lambda_k)$, $\underline{\lambda}'=
(\lambda'_1,\dots,\lambda'_k)$ be two $k$-tuples of disjoint properly
embedded arcs in $F$, with boundary in $\partial F\setminus Z$. 
We say that the pair $(\underline\lambda,\underline\lambda')$ is {\em
positive}, and write $\underline\lambda>\underline\lambda'$, if
along each component of $\partial F\setminus Z$ the points of 
$\partial(\bigcup_i \lambda'_i)$ all lie {\em before} those of
$\partial(\bigcup_i \lambda_i)$.
\end{defi}

\noindent
Similarly, given tuples 
$\underline\lambda^j=(\lambda^j_1,\dots,\lambda^j_k)$ ($j=0,\dots,\ell$),
we say that the sequence $(\underline\lambda^0,\dots,\underline\lambda^\ell)$
is positive if each pair $(\underline\lambda^j,\underline\lambda^{j+1})$ is positive,
i.e.\ $\underline\lambda^0>\dots>\underline\lambda^\ell$.

Given two tuples $\underline\lambda=(\lambda_1,\dots,\lambda_k)$ and
$\underline\lambda'=(\lambda'_1,\dots,\lambda'_k)$, we can perturb each arc
$\lambda_i$ (resp.\ $\lambda'_i$) by an isotopy that pushes it in the
positive (resp.\ negative) direction along $\partial F$, without crossing
$Z$, to obtain new tuples $\underline{\tilde\lambda}=(\tilde\lambda_1,
\dots,\tilde\lambda_k)$ and $\underline{\tilde{\lambda}}{}'=(
\tilde{\lambda}'_1,\dots,\tilde{\lambda}'_k)$ with the property that
$\underline{\tilde\lambda}>\underline{\tilde{\lambda}}{}'$.
Similarly, any sequence $(\underline\lambda^0,\dots,
\underline\lambda^\ell)$ can be made into a positive sequence by means of
suitable isotopies supported near $\partial F$ (again, the isotopies are
not allowed to cross $Z$).

\begin{figure}[t]
\setlength{\unitlength}{1cm}
\begin{picture}(6.05,3)(-1.05,-1.5)
\psset{unit=\unitlength}
\psellipticarc[linewidth=0.5pt,linestyle=dashed,dash=2pt 2pt](4.5,0)(0.2,1.5){90}{-90}
\psellipticarc(4.5,0)(0.2,1.5){-90}{90}
\psline[linearc=1.5](4.5,1.5)(-1,1.5)(-1,-1.5)(4.5,-1.5)
\pscircle*(4.58,1.3){0.07} \put(4.7,1.3){\small $z$}
\psellipse(0.35,0)(0.35,0.2)
\psellipse(1.65,0)(0.35,0.2)
\psellipticarc(3,0)(1.03,0.21){90}{180}
\psellipticarc(3,0)(1.71,0.51){90}{180}
\psellipticarc(3,0)(2.3,0.81){90}{180}
\psellipticarc(3,0)(3,1.11){90}{180}
\psline(3,0.2)(4.68,0.2)
\psline(3,0.5)(4.65,0.5)
\psline(3,0.8)(4.65,0.8)
\psline(3,1.1)(4.6,1.1)
\psellipticarc[linestyle=dashed,dash=2pt 2pt](3,0)(1.03,0.215){180}{270}
\psellipticarc[linestyle=dashed,dash=2pt 2pt](3,0)(1.71,0.515){180}{270}
\psellipticarc[linestyle=dashed,dash=2pt 2pt](3,0)(2.3,0.815){180}{270}
\psellipticarc[linestyle=dashed,dash=2pt 2pt](3,0)(3,1.115){180}{270}
\psline[linestyle=dashed,dash=2pt 2pt](3.05,-0.2)(4.32,-0.2)
\psline[linestyle=dashed,dash=2pt 2pt](3.05,-0.5)(4.32,-0.5)
\psline[linestyle=dashed,dash=2pt 2pt](3.05,-0.8)(4.35,-0.8)
\psline[linestyle=dashed,dash=2pt 2pt](3.05,-1.1)(4.38,-1.1)
\put(4.79,0.15){\small $\alpha_1$}
\put(4.75,1.05){\small $\alpha_{2g}$}
\psline[linestyle=dotted](4.9,0.85)(4.9,0.45)
\end{picture}
\qquad\quad
\begin{picture}(5.3,3)(-0.5,-1.5)
\psset{unit=\unitlength}
\psellipticarc[linewidth=0.5pt,linestyle=dashed,dash=2pt 2pt](4.5,0)(0.2,1.5){90}{-90}
\psellipticarc(4.5,0)(0.2,1.5){-90}{90}
\psline[linearc=1.5](4.5,1.5)(-1,1.5)(-1,-1.5)(4.5,-1.5)
\pscircle*(4.58,1.35){0.05} \put(4.7,1.3){\small $z$}
\psset{linecolor=red}
\psellipticarc[linestyle=dashed,dash=2pt 2pt](2.2,0)(0.23,0.215){180}{270}
\psellipticarc[linestyle=dashed,dash=2pt 2pt](2.2,0)(0.91,0.365){180}{270}
\psellipticarc[linestyle=dashed,dash=2pt 2pt](2.2,0)(1.5,0.515){180}{270}
\psellipticarc[linestyle=dashed,dash=2pt 2pt](2.2,0)(2.2,0.665){180}{270}
\pscurve[linestyle=dashed,dash=2pt 2pt](2.2,-0.2)(2.3,-0.15)(3.3,1.02)(4.35,1.17)
\pscurve[linestyle=dashed,dash=2pt 2pt](2.2,-0.35)(2.37,-0.3)(3.38,0.87)(4.35,1.02)
\pscurve[linestyle=dashed,dash=2pt 2pt](2.2,-0.5)(2.45,-0.45)(3.45,0.72)(4.3,0.87)
\pscurve[linestyle=dashed,dash=2pt 2pt](2.2,-0.65)(2.52,-0.6)(3.53,0.57)(4.3,0.72)
\pscurve[linestyle=dashed,dash=2pt 2pt](3.4,-1.5)(3.5,-1.45)(4.1,0.5)(4.3,0.57)
\pscurve[linestyle=dashed,dash=2pt 2pt](3.6,-1.5)(3.7,-1.45)(4.15,0.35)(4.3,0.42)
\pscurve[linestyle=dashed,dash=2pt 2pt](3.8,-1.5)(3.88,-1.45)(4.2,0.2)(4.3,0.27)
\pscurve[linestyle=dashed,dash=2pt 2pt](4,-1.5)(4.05,-1.45)(4.2,-0.1)(4.3,0.12)
\psset{linecolor=blue}
\psellipticarc(3,0)(1.03,0.81){90}{180}
\psellipticarc(3,0)(1.71,0.96){90}{180}
\psellipticarc(3,0)(2.3,1.11){90}{180}
\psellipticarc(3,0)(3,1.26){90}{180}
\psline(3,0.8)(4.65,0.8)
\psline(3,0.95)(4.63,0.95)
\psline(3,1.1)(4.6,1.1)
\psline(3,1.25)(4.58,1.25)
\psellipticarc[linestyle=dashed,dash=2pt 2pt](1.8,0)(1.8,1.5){180}{270}
\pscurve(1.8,-1.5)(2.3,-1.4)(3.6,0.5)(4.65,0.65)
\psellipticarc[linestyle=dashed,dash=2pt 2pt](2.1,0)(1.4,1.5){180}{270}
\pscurve(2.1,-1.5)(2.6,-1.4)(3.7,0.35)(4.68,0.5)
\psellipticarc[linestyle=dashed,dash=2pt 2pt](2.4,0)(1.1,1.5){180}{270}
\pscurve(2.4,-1.5)(2.9,-1.4)(3.8,0.2)(4.68,0.35)
\psellipticarc[linestyle=dashed,dash=2pt 2pt](2.7,0)(0.72,1.5){180}{270}
\pscurve(2.7,-1.5)(3.1,-1.4)(3.87,0.05)(4.7,0.2)
\psset{linecolor=red}
\psellipticarc(2.2,0)(0.23,0.215){90}{180}
\psellipticarc(2,0)(0.71,0.365){90}{180}
\psellipticarc(2.2,0)(1.5,0.515){90}{180}
\psellipticarc(2.2,0)(2.2,0.665){90}{180}
\pscurve(2.2,0.2)(2.3,0.18)(3.2,-1.45)(3.4,-1.5)
\pscurve(2,0.35)(2.4,0.3)(3.4,-1.45)(3.6,-1.5)
\pscurve(2.2,0.5)(2.6,0.42)(3.6,-1.45)(3.8,-1.5)
\pscurve(2.2,0.65)(2.7,0.55)(3.8,-1.45)(4,-1.5)
\psset{linecolor=black}
\psellipse(0.35,0)(0.35,0.2)
\psellipse(1.65,0)(0.35,0.2)
\pscircle*(0,0){0.04} \pscircle*(0.7,0){0.04} \pscircle*(1.3,0){0.04} \pscircle*(2,0){0.04}
\pscircle*(1.01,0.55){0.04} \pscircle*(1.75,0.63){0.04} \pscircle*(2.4,0.63){0.04} 
\pscircle*(1.5,0.44){0.04} \pscircle*(2.2,0.50){0.04} \pscircle*(2.09,0.35){0.04} 
\pscircle*(0.79,-0.5){0.04} \pscircle*(1.4,-0.61){0.04} \pscircle*(2.07,-0.65){0.04} 
\pscircle*(1.36,-0.41){0.04} \pscircle*(2.03,-0.5){0.04} \pscircle*(2.01,-0.35){0.04} 
\pscircle*(3.13,-0.07){0.04} \pscircle*(3.22,-0.33){0.04} 
\pscircle*(3.32,-0.66){0.04} \pscircle*(3.41,-0.89){0.04} 
\pscircle*(3.03,-0.26){0.04} \pscircle*(3.12,-0.59){0.04} 
\pscircle*(3.24,-0.92){0.04} \pscircle*(3.32,-1.11){0.04} 
\pscircle*(2.91,-0.56){0.04} \pscircle*(3.01,-0.87){0.04} 
\pscircle*(3.13,-1.13){0.04} \pscircle*(3.23,-1.27){0.04} 
\pscircle*(2.79,-0.78){0.04} \pscircle*(2.91,-1.07){0.04} 
\pscircle*(3.04,-1.28){0.04} \pscircle*(3.13,-1.38){0.04} 
\put(1.4,-1.8){\small $\tilde\alpha^0_{2g}\,\cdots\,\tilde\alpha^0_1$}
\put(3,-1.8){\small $\tilde\alpha^1_1\cdots\tilde\alpha^1_{2g}$}
\end{picture}
\caption{The arcs $\alpha_i$ and $\tilde\alpha_i^j$ on $(F,\{z\})$}
\label{fig:2g}
\end{figure}
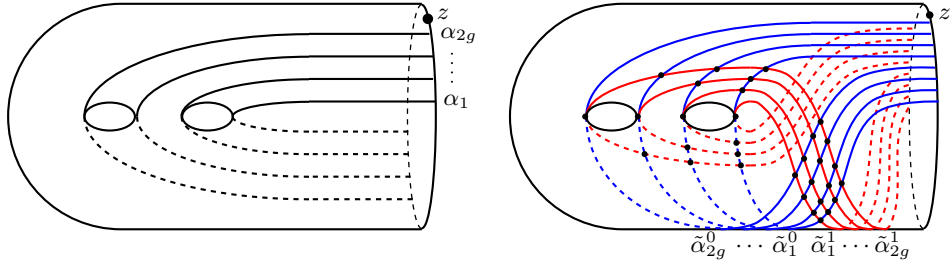

\begin{example*} Let $\underline\alpha=(\alpha_1,\dots,\alpha_{2g})$ be the
tuple of arcs represented on Figure~\ref{fig:2g} left: then the perturbed
tuples $\tilde{\underline\alpha}^j=(\tilde\alpha_1^j,
\dots,\tilde\alpha_{2g}^j)$ (Figure \ref{fig:2g} right) satisfy
$\tilde{\underline\alpha}^0>\tilde{\underline\alpha}^1$, i.e.\
the pair $(\tilde{\underline\alpha}^0,\tilde{\underline\alpha}^1)$ is
a positive perturbation of $(\underline\alpha,\underline\alpha)$.
\end{example*}

Next, consider a sequence $(L_0,\dots,L_\ell)$ of Lagrangian submanifolds
in the symmetric product $\Sym^k(F)$, each of which is either a closed
submanifold contained in the interior of $\Sym^k(F)$ or a product of
disjoint properly embedded arcs
$L_j=\lambda^j_1\times\dots\times\lambda^j_k$. Then we say that
the sequence $(L_0,\dots,L_\ell)$ is positive if, whenever $L_i$ and $L_j$
are products of disjointly embedded arcs for $i<j$, the corresponding
$k$-tuples of arcs satisfy
$\underline{\lambda}^i>\underline{\lambda}^j$. (There is no condition on the
closed Lagrangians).
 Modifying the arcs $\lambda^j_1,\dots,\lambda^j_k$
by suitable isotopies supported near $\partial F$ (without crossing 
$Z$) as above, given any sequence $(L_0,\dots,L_\ell)$ we can construct
Lagrangian submanifolds $\tilde{L}_0,\dots,\tilde{L}_\ell$ such that: (1)
$\tilde{L}_i$ is Hamiltonian isotopic to $L_i$, and either contained in the
interior of $\Sym^k(F)$ or a product of disjoint properly embedded arcs; and
(2) the sequence $(\tilde{L}_0,\dots,\tilde{L}_\ell)$ is positive.
We call $(\tilde{L}_0,\dots,\tilde{L}_\ell)$ a {\it positive perturbation}
of the sequence $(L_0,\dots,L_\ell)$.

With this understood, we can now give an informal (and imprecise) definition of the 
partially wrapped Fukaya category of the symmetric product $\Sym^k(F)$ 
relative to the set $Z$; we are still assuming that every component of
$\partial F$ contains at least one point of $Z$. 
The reader is referred to \cite{fuksymg} for a more precise construction.

\begin{defi}
The partially wrapped Fukaya category $\F=\F(\Sym^k(F),Z)$ is an
$A_\infty$-category with objects of
two types:
\begin{enumerate}
\item closed balanced Lagrangian submanifolds lying in the interior
of $\Sym^k(F)$;
\item properly embedded Lagrangian submanifolds of the form $\lambda_1
\times\dots\times\lambda_k$, where $\lambda_i$ are disjoint properly
embedded arcs with boundary contained in $\partial F\setminus Z$.
\end{enumerate}
Morphism spaces and compositions are defined by perturbing objects
of the second type in a suitable manner near the boundary so that they 
form positive sequences. Namely, we set
$\hom_\F(L_0,L_1)=CF(\tilde{L}_0,\tilde{L}_1)$ (i.e., the $\Z_2$-vector space generated by
points of $\tilde{L}_0\cap \tilde{L}_1$, with a differential counting rigid
holomorphic discs), where
$(\tilde{L}_0,\tilde{L}_1)$ is a suitably chosen positive perturbation of the pair $(L_0,L_1)$.
The composition $m_2:\hom_\F(L_0,L_1)\otimes\hom_\F(L_1,L_2)\to
\hom_\F(L_0,L_2)$
 and higher products $m_\ell:\hom_\F(L_0,L_1)\otimes\dots\otimes
\hom_\F(L_{\ell-1},L_\ell)\to \hom_\F(L_0,L_\ell)$ are similarly defined
by perturbing $(L_0,\dots,L_\ell)$ to a positive sequence $(\tilde{L}_0,
\dots,\tilde{L}_\ell)$ and counting rigid holomorphic discs with boundary 
on the perturbed Lagrangians.
\end{defi} 

The extended category $\F^\#=\F^\#(\Sym^k(F),Z)$ is defined
similarly, but also includes closed balanced generalized Lagrangian
submanifolds of $\Sym^k(F)$ (i.e., formal images of Lagrangians under
sequences of balanced Lagrangian correspondences) of the sort introduced in
\S \ref{s:LP}.
\medskip

To be more precise, the construction of the partially wrapped Fukaya
category involves the completion $\hat{F}=F\cup (\partial F\times
[1,\infty))$, and its symmetric product $\Sym^k(\hat{F})$. 
Arcs in $F$ can be completed to properly embedded arcs in $\hat{F}$,
translation-invariant in the cylindrical ends, and hence the objects
of $\F(\Sym^k(F),Z)$ can be completed to properly embedded Lagrangian
submanifolds of $\Sym^k(\hat{F})$ which are cylindrical at infinity.
The Riemann surface $\hat{F}$ carries a Hamiltonian vector field supported 
away from the interior of $F$ and whose positive (resp.\ negative) time
flow rotates the cylindrical ends of $\hat{F}$ in the positive (resp.\
negative) direction and accumulates towards the rays $Z\times [1,\infty)$.
(In the cylindrical ends $\partial F\times [1,\infty)$, the generating Hamiltonian function
$h$ is of the form $h(x,r)=\rho(x)r$ where $\rho:\partial F\to [0,1]$ satisfies
$\rho^{-1}(0)=Z$). This flow on $\hat{F}$ can be used to construct a 
Hamiltonian flow on $\Sym^k(\hat{F})$ which preserves the product 
structure away from the diagonal (namely, the generating Hamiltonian
is given by $H(\{z_1,\dots,z_k\})=\sum_i h(z_i)$ away from the diagonal).
The $A_\infty$-category $\F(\Sym^k(F),Z)$ is then constructed in essentially
the same manner as the wrapped Fukaya category defined by Abouzaid and
Seidel~\cite{AS}: namely, morphism spaces are limits of the Floer complexes
upon long-time perturbation by the Hamiltonian flow.
(In general various technical issues could arise with this construction, but
product Lagrangians in $\Sym^k(\hat{F})$ are fairly well-behaved, see 
\cite{fuksymg}).
\medskip

When a component of $\partial F$ does not contain any point of $Z$, 
the Hamiltonian flow that we consider rotates the 
corresponding cylindrical end of $\hat{F}$ by arbitrarily large amounts. Hence
the perturbation causes properly embedded arcs in $\hat{F}$
to wrap around the cylindrical end infinitely many times, which typically
makes the complex $\hom_\F(L_0,L_1)$ infinitely generated when
$L_0$ and $L_1$ are non-compact objects of $\F(\Sym^k(F),Z)$. For instance, when $Z=\emptyset$ the
category we consider is simply the wrapped Fukaya category of
$\Sym^k(\hat{F})$ as defined in \cite{AS}.

\subsection{The algebra of a decorated surface} 

\begin{defi}
A {\em decorated surface} is a triple $\FF=(F,Z,\underline{\alpha})$ where
$F$ is a connected compact Riemann surface with non-empty boundary,
$Z$ is a finite subset of $\partial F$, and $\underline{\alpha}=\{\alpha_1,
\dots,\alpha_n\}$ is a collection of disjoint properly embedded arcs in $F$
with boundary in $\partial F\setminus Z$.
\end{defi}

Given a decorated surface $\FF=(F,Z,\underline\alpha)$, an integer 
$k\le n$, and a $k$-element subset $s\subseteq \{1,\dots,n\}$, 
the product $D_s=\prod_{i\in s}\alpha_i$ is an object of 
$\F=\F(\Sym^k(F),Z)$. The endomorphism algebra of the direct sum of these
objects is an $A_\infty$-algebra naturally associated to $\FF$.

\begin{defi}
For $k\le n$, denote by $\SS^n_k$ the set of all $k$-element subsets of 
$\{1,\dots,n\}$. Then to a decorated surface $\FF=(F,\,Z,\,\underline\alpha=
\{\alpha_1,\dots,\alpha_n\})$
and an integer $k\le n$ we associate the $A_\infty$-algebra
$$\textstyle \A(\FF,k)=\bigoplus\limits_{s,t\in\SS^n_k}\hom_{\F}(D_s,D_t),\quad
\text{where}\ D_s=\prod\limits_{i\in s}\alpha_i,$$
with differential and products defined by those of $\F=\F(\Sym^k(F),Z)$.
\end{defi}

In the rest of this section, we focus on a special case where $\A(\FF,k)$
can be expressed in purely combinatorial terms, and is in fact isomorphic 
to (the obvious generalization of) the
{\it bordered algebra} introduced by Lipshitz, Ozsv\'ath and 
Thurston~\cite{LOT}. The following proposition implies 
Theorem \ref{thm:algebra}
as a special case:

\begin{proposition}\label{prop:nice}
Let $\FF=(F,Z,\underline\alpha)$ be a decorated surface, and assume that 
every connected component
of $F\setminus (\alpha_1\cup\dots\cup\alpha_n)$ contains at least one
point of $Z$. For $i,j\in\{1,\dots,n\}$, denote by $\chi_i^j$ the
set of {\em chords} from $\partial\alpha_i$ to $\partial\alpha_j$ in
$\partial F\setminus Z$, i.e.\ homotopy classes of immersed arcs
$\gamma:[0,1]\to \partial F\setminus Z$ such that $\gamma(0)\in\partial
\alpha_i$, $\gamma(1)\in\partial\alpha_j$, and the tangent vector $\gamma'(t)$ is always
oriented in the positive direction along $\partial F$. Moreover, denote
$\bar\chi{}_i^i$ the set obtained by adjoining to $\chi_i^i$ an extra
element~$\mathbf{1}_i$, and let $\bar\chi{}_i^j=\chi_i^j$ for $i\neq j$.
Then the following properties hold:
\begin{itemize}
\item Given $s,t\in \SS^n_k$, let $s=\{i_1,\dots,i_k\}$, and
denote by $\Phi(s,t)$ the set of bijective maps from $s$ to $t$.
Then the $\Z_2$-vector space $\hom_{\F(\Sym^k(F),Z)}(D_s,D_t)$
admits a basis indexed by the elements of $$\bar\chi_s^t:=
\smash{\bigsqcup\limits_{f\in \Phi(s,t)}}
\left(\bar\chi_{i_1}^{f(i_1)}
\times\dots\times \bar\chi_{i_k}^{f(i_k)}\right).$$
\item The differential and product in $\A(\FF,k)$ are determined
by explicit combinatorial formulas as in \cite{LOT}.
\item The higher products $\{m_\ell\}_{\ell\ge 3}$ vanish identically,
i.e.\ the $A_\infty$-algebra $\A(\FF,k)$ is in fact a differential algebra.
\end{itemize}
\end{proposition}

\begin{figure}[t]
\setlength{\unitlength}{9mm}
\hfil \begin{picture}(5,4.8)(0,-2.4)
\psset{unit=\unitlength}
\multips(0,-0.6)(0,0.4){4}{\psline(0,0)(4,0)}
\multips(0,-0.6)(0,0.4){4}{\pscurve[linestyle=dashed,dash=8pt 2pt](0,0.06)(0.8,0.08)(1,0.11)(3,1.55)(3.2,1.57)(4,1.6)}
\multips(0,-0.6)(0,0.4){4}{\pscurve[linestyle=dashed,dash=2pt 2pt](0,-0.06)(0.8,-0.08)(1,-0.11)(3,-1.55)(3.2,-1.57)(4,-1.6)}
\psline(4,-2.5)(4,2.5)
\put(4.2,1.6){\makebox(0,0)[lc]{$\Biggr\}\ \underline{\tilde\alpha}^0$}}
\put(4.2,0){\makebox(0,0)[lc]{$\Biggr\}\ \underline{\tilde\alpha}^1$}}
\put(4.2,-1.6){\makebox(0,0)[lc]{$\Biggr\}\ \underline{\tilde\alpha}^2$}}
\pscircle*(0.2,0.4){0.05}
\pscircle*(0.2,0){0.05}
\pscircle*(0.2,-0.4){0.05}
\pscircle*(0.2,0.8){0.05}
\pscircle*(0.2,-0.8){0.05}
\end{picture}
\caption{Positive perturbations of $\underline{\alpha}$ near $\partial F$}
\label{fig:nice}
\end{figure}
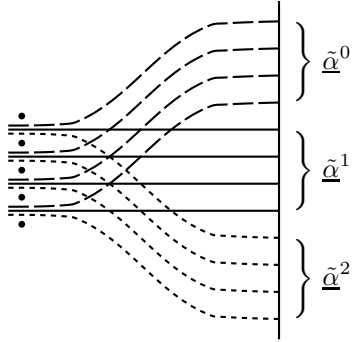

\proof[Sketch of proof\/ \rm (see also \cite{fuksymg}).]
For $\ell\ge 1$, we construct perturbations 
$\tilde{\underline\alpha}^0,\dots,\tilde{\underline\alpha}^\ell$
of $\underline\alpha$, with
$\underline{\tilde\alpha}^0>\dots>\underline{\tilde\alpha}^\ell$,
in such a way that the diagram formed by the $\ell+1$ collections of
$n$ arcs $\tilde\alpha^j_i$ on $F$ enjoys properties similar to those of
``nice'' diagrams in Heegaard-Floer theory (cf.\ \cite{SW}). Namely, we ask that
for each $i$ the arcs $\tilde\alpha^0_i,\dots,
\tilde\alpha^\ell_i$ remain close to $\alpha_i$ in the interior of
$F$, where any two of them intersect transversely exactly once; the
total number of intersections in the diagram is minimal; and
all intersections between the arcs of $\underline{\tilde\alpha}^j$ and
those of $\underline{\tilde\alpha}^{j'}$ are transverse and occur with the same
oriented angle 
$(j-j')\theta$ (for a fixed small $\theta>0$) between the two arcs
at the intersection point. Hence the local picture near any interval
component of $\partial F\setminus Z$ is as shown in Figure \ref{fig:nice}.
(At a component of $\partial F$ which does not carry a point of $Z$, we
need to consider arcs which wrap
infinitely many times around the cylindrical end of the completed surface
$\hat{F}$, but the situation is otherwise unchanged).

For $j<j'$ and $i,i'\in \{1,\dots,n\}$ we
have a natural bijection between the points of 
$\tilde\alpha^j_i\cap \tilde\alpha^{j'}_{i'}$ and the elements of
$\bar\chi_i^{i'}$. Hence, passing to the symmetric product, 
the intersections of $\tilde{D}_s^j=\prod_{i\in s}
\tilde\alpha_i^j$ and $\tilde{D}_t^{j'}=\prod_{i'\in t}
\tilde\alpha_{i'}^{j'}$ are transverse and
in bijection with the elements of $\bar\chi_s^t$. The first claim follows.

The rest of the proposition follows from a calculation of the Maslov index of a
holomorphic disc in $\Sym^k(F)$ with boundary on $\ell+1$ product
Lagrangians $\tilde{D}^0_{s_0},\dots,\tilde{D}^\ell_{s_\ell}$. Namely,
let $\phi$ be the homotopy class of such a holomorphic disc contributing to 
the order $\ell$ product in $\A(\FF,k)$. Projecting from the symmetric
product to $F$, we can think of $\phi$ as a 2-chain in $F$ with boundary on
the arcs of the diagram, staying within the bounded regions of the diagram 
(i.e., those which do not intersect $\partial F$). Then the Maslov index
$\mu(\phi)$ and the intersection number $i(\phi)$ of $\phi$ with the
diagonal divisor in $\Sym^k(F)$ are related to each other by the following
formula due to Sarkar \cite{Sarkar}:
\begin{equation}\label{eq:sarkar}
\mu(\phi)=i(\phi)+2e(\phi)-(\ell-1)k/2,
\end{equation}
where $e(\phi)$ is the {\it Euler measure} of the 2-chain $\phi$,
characterized by additivity and by the property that the Euler 
measure of an embedded $m$-gon with convex corners is $1-\frac{m}{4}$.
On the other hand, since every component of $F\setminus
(\alpha_1\cup\dots\cup\alpha_n)$ contains a point of $Z$,
the regions of the diagram corresponding to those components remain
unbounded after perturbation. In particular, the regions marked by dots
in Figure~\ref{fig:nice} are all unbounded, and hence not part of the
support of $\phi$. 

This implies that the support of $\phi$ is contained
in a union of regions which are either planar (as in Figure \ref{fig:nice}) 
or cylindrical (in the case of a component of $\partial F$ which does not
carry any point of $Z$), and within which the Euler measure of a convex
polygonal region can be computed by summing contributions from its vertices,
namely $\frac{1}{4}-\frac{\vartheta}{2\pi}$ for a vertex with angle $\vartheta$.
Considering the respective contributions of the $(\ell+1)k$ corners of
the chain $\phi$ (and observing that the contributions from any other 
vertices traversed by the boundary of $\phi$ cancel out), we conclude that
$e(\phi)=(\ell-1)k/4$, and $\mu(\phi)=i(\phi)\ge 0$.

On the other hand, $m_\ell$ counts rigid holomorphic discs, for which
$\mu(\phi)=2-\ell$. This immediately implies that $m_\ell=0$
for $\ell\ge 3$. For $\ell=1$, the diagram we consider is ``nice'', i.e.\
all the bounded regions are quadrilaterals;
as observed by Sarkar and Wang, this implies that the Floer differential
on $CF(\tilde{D}_s^0,\tilde{D}_t^1)$ counts empty embedded rectangles
\cite[Theorems~3.3 and~3.4]{SW}. Finally, for $\ell=2$, the Maslov index
formula shows that the product counts discs which are disjoint
from the diagonal strata in $\Sym^k(F)$. By an argument similar to that
in \cite{LMW} (see also \hbox{\cite[Proposition 3.5]{fuksymg}}), this implies
that $m_2$ counts $k$-tuples of immersed holomorphic triangles in $F$ which 
either are disjoint or overlap head-to-tail (cf.\ \cite[Lemma~2.6]{LMW}).

Finally, these combinatorial descriptions of $m_1$ and $m_2$ in terms
of diagrams on $F$ can be recast in terms of Lipshitz, Ozsv\'ath and
Thurston's definition of differentials and products in the bordered algebra
\cite{LOT}. Namely, the dictionary between points of $\tilde{D}_s^0\cap
\tilde{D}_t^1$ proceeds by matching intersections of $\tilde\alpha_i^0$ with
$\tilde\alpha_{j}^1$ near $\partial F$ with chords from $\alpha_i$ to
$\alpha_j$ (pictured as upwards strands in the notation of \cite{LOT}), and the 
intersection of $\tilde\alpha_i^0$ with $\tilde\alpha_i^1$ in the interior of 
$F$ with a pair of horizontal dotted lines in the graphical notation of
\cite{LOT}. See \cite[section 3]{fuksymg} for details.
\endproof

\subsection{Generating the partially wrapped Fukaya category}
\label{ss:generate}

In this section, we outline the proof of Theorem
\ref{thm:generate}. The main ingredients are Lefschetz fibrations on the symmetric
product, their Fukaya categories as defined and studied by Seidel
\cite{SeVCM,SeBook}, and acceleration 
$A_\infty$-functors between partially wrapped Fukaya 
categories.

\subsubsection{Lefschetz fibrations on the symmetric product}

Let $\hat{F}$ be an open Riemann surface (with infinite cylindrical ends),
and let $\pi:\hat{F}\to\C$ be
a branched covering map. Assume that the critical points
$q_1,\dots,q_n$ of $\pi$ are non-degenerate (i.e., the covering $\pi$ is
{\it simple}), and that
the critical values $p_1,\dots,p_n\in\C$ are distinct, lie in the unit disc, 
and satisfy $\Im(p_1)<\dots<\Im(p_n)$.

Each critical point $q_j$ of $\pi$ determines a properly embedded arc
$\hat{\alpha}_j\subset \hat{F}$, namely the union of the two lifts of the
half-line $\R_{\ge 0}+p_j$ which pass through $q_j$. 

We consider the $k$-fold symmetric product $\Sym^k(\hat{F})$
($1\le k\le n$), equipped with the product complex structure $J$,
and the
holomorphic map $f_{n,k}:\Sym^k(\hat{F})\to \C$ defined by
$f_{n,k}(\{z_1,\dots,z_k\})=\pi(z_1)+\dots+\pi(z_k)$.

\begin{proposition}\label{prop:lf}
$f_{n,k}:\Sym^k(\hat{F})\to \C$ is a holomorphic map with isolated
non-degenerate critical points (i.e., a Lefschetz fibration);
its $\binom{n}{k}$ critical points are the tuples consisting
of $k$ distinct points in $\{q_1,\dots,q_{n}\}$.
\end{proposition}

\proof
Given $\underline{z}\in\Sym^k(\hat{F})$, denote by
$z_1,\dots,z_r$ the distinct elements in the $k$-tuple
$\underline{z}$, and by $k_1,\dots,k_r$ the multiplicities with which
they appear. The tangent space $T_{\underline{z}}\Sym^k(\hat{F})$ decomposes
into the direct sum of the $T_{\{z_i,\dots,z_i\}}\Sym^{k_i}(\hat{F})$,
and $df_{n,k}(\underline{z})$ splits into the direct sum of the
differentials $df_{n,k_i}(\{z_i,\dots,z_i\})$. Thus $\underline{z}$ is a
critical point of $f_{n,k}$ if and only if $\{z_i,\dots,z_i\}$ is a critical
point of $f_{n,k_i}$ for each $i\in \{1,\dots,r\}$.

By considering the 
restriction of $f_{n,k_i}$ to the diagonal stratum, we see that
$\{z_i,\dots,z_i\}$ cannot be a critical point of $f_{n,k_i}$ unless $z_i$
is a critical point of $\pi$. Assume now that $z_i$ is a critical point
of $\pi$, and pick a local complex coordinate $w$ on $\hat{F}$ near $z_i$,
in which $\pi(w)=w^2+\mathrm{constant}$. Then a neighborhood
of $\{z_i,\dots,z_i\}$ in $\Sym^{k_i}(\hat{F})$ identifies with a 
neighborhood of the origin in $\Sym^{k_i}(\C)$, with coordinates
given by the elementary symmetric functions $\sigma_1,\dots,\sigma_{k_i}$.
The local model for $f_{n,k_i}$
is then
$$f_{n,k_i}(\{w_1,\dots,w_{k_i}\})=w_1^2+\dots+w_{k_i}^2+\mathrm{constant}=
\sigma_1^2-2\sigma_2+\mathrm{constant}.$$ Thus, for $k_i\ge 2$ the point
$\{z_i,\dots,z_i\}$
is never a critical point of $f_{n,k_i}$. We conclude that the only critical
points of $f_{n,k}$ are tuples of distinct critical points of $\pi$;
moreover
these critical points are clearly non-degenerate.
\endproof

For $s\in \mathcal{S}^n_k$, we denote by $Q_s=\{q_i,\ i\in s\}$ the
corresponding critical point of $f_{n,k}$, and by $P_s=\sum_{i\in s}
p_i$ the associated critical value.

As in \S \ref{s:LP}, equip
$\Sym^k(\hat{F})$ with a K\"ahler form $\omega$ which is of product type 
away from the diagonal strata, and the associated K\"ahler metric. 
This allows us to associate
to each critical point $Q_s$ a properly embedded Lagrangian disc $\hat{D}_s$
in $\Sym^k(\hat{F})$ (called {\it Lefschetz thimble}), namely the
set of those points in \hbox{$f_{n,k}^{-1}(\R_{\ge 0}+P_s)$} for which the
gradient flow of $\Re f_{n,k}$ converges to the critical point $Q_s$.
A straightforward calculation shows that $\hat{D}_s=\prod_{i\in s}
\hat{\alpha}_i$.

More generally, one can associate a Lefschetz thimble to any properly
embedded arc $\gamma$ connecting $P_s$ to infinity: namely, the set of points
in $f_{n,k}^{-1}(\gamma)$ for which symplectic parallel transport converges 
to the critical point $Q_s$. We will only consider the case where 
$\gamma$ is a straight half-line. Given $\theta\in
(-\frac{\pi}{2},\frac{\pi}{2})$, the thimble
associated to the half-line $e^{i\theta}\R_{\ge 0}+P_s$ is again a product
$\hat{D}_s(\theta)=\prod_{i\in s} \hat{\alpha}_i(\theta)$, where $\hat\alpha_i(\theta)$ is
the
union of the two lifts of the half-line $e^{i\theta}\R_{\ge 0}+p_j$ through
$q_j$.

\subsubsection{A special case of Theorem \ref{thm:generate}}\label{ss:332}
In the same setting as above, consider the Riemann surface with
boundary $F=\pi^{-1}(D^2)$, i.e.\ the preimage of the unit disc, and
let $Z=\pi^{-1}(-1)\subset \partial F$. Let 
$\alpha_i=\hat\alpha_i\cap F$, and $D_s=\hat{D}_s\cap \Sym^k(F)=
\prod_{i\in s}\alpha_i$.
Then we can reinterpret the partially wrapped Fukaya category
$\F(\Sym^k(F),Z)$ and the algebra $\A(\FF,k)$ associated to the arcs $\alpha_1,\dots,
\alpha_n$ in different terms.

Seidel associates to the Lefschetz fibration $f_{n,k}$ a Fukaya category
$\F(f_{n,k})$, whose objects are compact Lagrangian submanifolds of
$\Sym^k(\hat{F})$ on one hand, and Lefschetz thimbles associated to
admissible arcs connecting a critical value of $f_{n,k}$ to infinity on 
the other hand \cite{SeBook}. Here we say that an arc is admissible
with slope $\theta\in(-\frac{\pi}{2},\frac{\pi}{2})$ if outside of a compact
set it is a half-line of slope $\theta$. (Seidel
considers the case of an exact symplectic form, and defines things somewhat
differently; however our setting
does not pose any significant additional difficulties). 

Morphisms between
thimbles in $\F(f_{n,k})$ (and compositions thereof) are defined by
means of suitable perturbations. Namely, given two admissible arcs $\gamma_0,
\gamma_1$ and the corresponding thimbles $D_0,D_1\subset\Sym^k(\hat{F})$, 
one sets $\hom_{\F(f_{n,k})}(D_0,D_1)=CF(\tilde{D}_0,\tilde{D}_1)$, where
$\tilde{D}_0,\tilde{D}_1$ are thimbles obtained by suitably perturbing
$(\gamma_0,\gamma_1)$ to a positive pair $(\tilde\gamma_0,
\tilde\gamma_1)$, i.e.\ one for which the slopes satisfy
$\theta'_0>\theta'_1$.

Restricting ourselves to the special case of straight half-lines,
and observing that for sufficiently small $\theta_0>\dots>\theta_\ell$ the
collections of arcs $\alpha_i(\theta_j)=\hat\alpha_i(\theta_j)\cap F$ form
a positive sequence in the sense of \S \ref{ss:wrap}, it is not hard to see
that we have an isomorphism of $A_\infty$-algebras
$$\textstyle\bigoplus\limits_{s,t\in \SS^n_k}\hom_{\F(\Sym^k(F),Z)}(D_s,D_t)\simeq
\bigoplus\limits_{s,t\in \SS^n_k}\hom_{\F(f_{n,k})}(\hat{D}_s,\hat{D}_t).$$

A key result due to Seidel is the following:

\begin{theorem}[Seidel \cite{SeBook}, Theorem 18.24]
The $A_\infty$-category $\F(f_{n,k})$ is generated by the exceptional
collection of thimbles $\hat{D}_s$, $s\in\SS^n_k$.
\end{theorem}

\noindent In other terms, every object of $\F(f_{n,k})$ is quasi-isomorphic
to a twisted complex built out of the objects $\hat{D}_s$, $s\in \SS^n_k$.

This implies Theorem \ref{thm:generate} in the special case where
$F$ is a simple branched cover of the disc with $n$ critical points, 
$Z$ is the preimage of $-1$, and the arcs $\alpha_1,\dots,\alpha_n$ 
are lifts of half-lines connecting 
connecting the critical values to
the boundary of the disc along the real positive direction. 
(More precisely, in view of the relation
between $\F(f_{n,k})$ and $\F(\Sym^k(F),Z)$, Seidel's result directly 
implies that the compact objects of $\F(\Sym^k(F),Z)$ are generated by 
the $D_s$. On the other hand, arbitrary products of properly embedded
arcs cannot be viewed as objects of $\F(f_{n,k})$, but by performing
sequences of arc slides we can express them explicitly as iterated
mapping cones involving the generators $D_s$, see below.)

\subsubsection{Acceleration functors}\label{ss:333}
Consider a fixed surface $F$, and two subsets $Z\subseteq Z'\subset \partial
F$. Then there exists a natural $A_\infty$-functor from $\F(\Sym^k(F),Z')$
to $\F(\Sym^k(F),Z)$, called ``acceleration functor''. This functor is
identity on objects, and in the present case it is simply given by an inclusion 
of morphism spaces.
In general, it is given by the Floer-theoretic continuation maps that arise
when comparing the Hamiltonian perturbations used to define morphisms and
compositions in $\F(\Sym^k(F),Z')$ and $\F(\Sym^k(F),Z)$.

\begin{figure}[t]
\setlength{\unitlength}{9mm}
\hfil \begin{picture}(5,3)(0,-1.5)
\psset{unit=\unitlength}
\psline(0,-0.6)(4,-0.6)
\psline(0,0.6)(4,0.6)
\pscircle*(4,0){0.05}
\put(4.2,0.1){\small $z\in Z'\setminus Z$}
\put(0.5,-0.42){$\delta_i$}
\put(0.5,0.85){$\lambda_j$}
\psline(4,-1.5)(4,1.5)
\put(6.5,0){$\to$}
\end{picture}
\qquad\qquad\qquad
\begin{picture}(5,3)(0,-1.5)
\psset{unit=\unitlength}
\psline(4,-1.5)(4,1.5)  
\pscurve(0,0.54)(0.8,0.52)(1,0.49)(3,-0.95)(3.2,-0.97)(4,-1)
\pscurve(0,-0.54)(0.8,-0.52)(1,-0.49)(3,0.95)(3.2,0.97)(4,1)
\put(0.5,-0.36){$\tilde\delta_i$}
\put(0.5,0.8){$\tilde\lambda_j$}
\put(1.65,0.3){\small I}
\put(1.6,-0.5){\small II}
\pscircle*(1.7,0){0.05}
\end{picture}
\caption{Perturbation and the acceleration functor}
\label{fig:accelerate}
\end{figure}
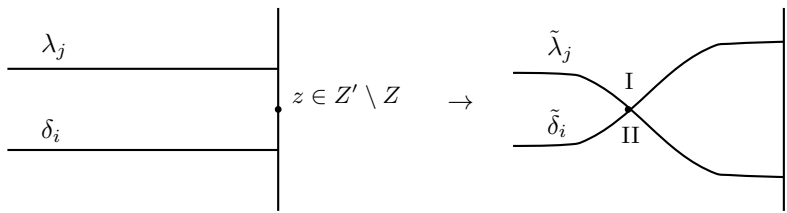

Consider two products $\Delta=\delta_1\times\dots\times\delta_k$ and
$L=\lambda_1\times\dots\times\lambda_k$ of disjoint properly embedded
arcs in $F$ with boundary in $\partial F\setminus Z'$. Perturbing the arcs
$\delta_1,\dots,\delta_k$ and $\lambda_1,\dots,\lambda_k$ near $\partial F$
if needed (without crossing $Z'$), we can assume that the pair $(\Delta,L)$ 
is positive with respect to $Z'$. On the other hand, achieving positivity
with respect to the smaller subset $Z$ may require a further perturbation
of the arcs $\delta_i$ (resp.\ $\lambda_i$) in the positive (resp.\
negative) direction along $\partial F$, to obtain product Lagrangians 
$\tilde\Delta=\tilde{\delta}_1\times\dots\times\tilde\delta_k$ and
$\tilde{L}=\tilde\lambda_1\times\dots\times\tilde\lambda_k$. This
perturbation can be performed in such a way as to only {\it create}
new intersection points. The local picture is as shown on Figure
\ref{fig:accelerate}. The key observation is that none of the intersection
points created in the isotopy can be the outgoing end of a holomorphic
strip in $\Sym^k(F)$ with boundary on $\tilde{\Delta}\cup\tilde{L}$ and
whose incoming end is a previously existing intersection point (i.e.,
one that arises by deforming a point of $\Delta\cap L$). Indeed, considering
Figure \ref{fig:accelerate} right, locally the projection of this 
holomorphic strip to $F$ would cover one of the two regions labelled I and
II; but then by the maximum principle it would need to hit $\partial F$,
which is not allowed. This implies that $CF(\Delta,L)$ is
naturally a subcomplex of $CF(\tilde\Delta,\tilde{L})$. The same argument
also holds for products and higher compositions, ensuring that the
acceleration functor is well-defined.

In particular, given a collection $\underline{\alpha}=(\alpha_1,\dots,\alpha_n)$
of disjoint properly embedded arcs  in $F$, and setting
$\FF=(F,Z,\underline\alpha)$ and $\FF'=(F,Z',\underline\alpha')$, we
obtain that $\A(\FF',k)$ is naturally an $A_\infty$-subalgebra of
$\A(\FF,k)$ for all $k$.

Finally, one easily checks that the acceleration functor is unital (at
least on cohomology), and surjective on (isomorphism classes of) 
objects. Hence, if the $\binom{n}{k}$ objects $D_s=\prod_{i\in s} \alpha_i$
$(s\in \SS^n_k)$
generate $\F(\Sym^k(F),Z')$, then they also generate $\F(\Sym^k(F),Z)$.
(Indeed, the assumption means that any object $L$ of $\F(\Sym^k(F),Z')$ is 
quasi-isomorphic to a twisted complex built out of the $D_s$; since
$A_\infty$-functors are exact, this implies that $L$ is also 
quasi-isomorphic to the corresponding twisted complex in
$\F(\Sym^k(F),Z)$).

\subsubsection{Eliminating generators by arc slides}
We now consider a general decorated surface $\FF=(F,Z,\underline\alpha)$.
The arcs $\alpha_1,\dots,\alpha_n$ on $F$ might not be a full set of Lefschetz
thimbles for any simple branched covering map, but they are always a 
subset of the thimbles
of a more complicated covering (with $m$ critical points, $m\geq n$).
Namely, after a suitable deformation (which does
not affect the symplectic topology of the completed symmetric product
$\Sym^k(\hat{F})$),
we can always assume that $F$ projects to the disc by a simple branched 
covering map $\pi$ with critical values $p_1,\dots,p_m$, in such a way 
that the arcs $\alpha_1,\dots,\alpha_n$
are lifts of $n$ of the half-lines $\R_{\ge 0}+p_j$,
while each point
of $Z$ projects to $-1$. Hence, taking the remaining critical
values of $\pi$ and elements of $\pi^{-1}(-1)$ into account, there
exists a subset $Z'\supseteq Z$ of $\partial F$, and
a collection $\underline\alpha'$ of $m\ge n$ disjoint properly embedded arcs (including the
$\alpha_i$), such that $\FF'=(F,Z',\underline\alpha')$ is as in \S \ref{ss:332}. 
Then, as seen above, the partially wrapped Fukaya category
$\F(\Sym^k(F),Z')$ is generated by the $\binom{m}{k}$ product objects
$D'_s=\prod_{i\in s}\alpha'_i$ ($s\in \SS^m_k$). 

Moreover, by considering
the acceleration functor as in \S \ref{ss:333}, we conclude that
$\F(\Sym^k(F),Z)$ is also generated by the objects $D'_s$, $s\in \SS^m_k$.
Thus, Theorem \ref{thm:generate} follows if, assuming
that each component of $F\setminus (\alpha_1\cup\dots\cup \alpha_n)$
is a disc containing at most one point of $Z$, we can show that the
$\binom{m}{k}-\binom{n}{k}$ additional objects we have introduced can be
expressed in terms of the others. This is done by eliminating the
additional arcs $\alpha'_i$ one at a time.

\begin{figure}[t]
\setlength{\unitlength}{9mm}
\hfil
\begin{picture}(12,4)(0,0)
\psset{unit=\unitlength}
\psellipse[linewidth=0.5pt,linestyle=dashed,dash=2pt 2pt](5.5,2)(0.2,2)
\psellipticarc(5.5,2)(0.2,2){-90}{90}
\psellipticarc(4,2.75)(0.4,0.3){40}{140}
\psellipticarc(4,3.25)(0.6,0.4){-140}{-40}
\psellipticarc(2,1.8)(0.4,0.3){40}{140}
\psellipticarc(2,2.3)(0.6,0.4){-140}{-40}
\psellipticarc(3.5,0.45)(0.4,0.3){40}{140}
\psellipticarc(3.5,0.95)(0.6,0.4){-140}{-40}
\psline(5.5,4)(2.5,4) \psline(5.5,0)(2.5,0)
\psellipticarc(2.5,2)(2.5,2){90}{-90}
\psellipticarc(5.6,3)(3,0.4){90}{-90}
\psellipticarc(5.6,0.6)(3.5,0.4){90}{-90}
\psellipticarc[linestyle=dashed,dash=4pt 2pt](5.6,0.75)(3.5,0.4){90}{-90}
\psellipticarc(3.6,2)(2.8,0.4){90}{-90}
\psline(3.6,2.4)(5.7,2.4)
\psline(3.6,1.6)(5.7,1.6)
\psellipticarc[linestyle=dashed,dash=4pt 2pt](3.6,2.15)(2.8,0.45){120}{-90}
\psline[linestyle=dashed,dash=4pt 2pt](2.2,2.52)(5.6,2.75)
\psline[linestyle=dashed,dash=4pt 2pt](3.6,1.7)(5.6,1.7)
\psellipticarc(5.6,3)(3.6,0.6){90}{180}
\psellipticarc(1,3)(1,0.4){-70}{0}
\psellipticarc(1.8,2.1)(1.3,0.6){110}{180}
\psellipticarc(3.6,2.1)(3.1,0.7){180}{270}
\psline(3.6,1.4)(5.7,1.4)
\put(1.4,3){\small $\lambda''_1$}
\put(2.2,2.8){\small $\lambda'_1$}
\put(2.6,2.08){\small $\lambda_1$}
\put(3.6,1.88){\small $\tilde\lambda_1$}
\put(1.75,0.98){\small $\tilde\lambda_2$}
\put(2.1,0.2){\small $\lambda_2$}
\pscircle*(2.2,0.68){0.07}
\put(1.6,0.6){\small $x_2$}
\pscircle*(4.35,2.65){0.07}
\put(4.5,2.82){\small $x'_1$}
\psline{->}(6.5,2)(7.5,2)
\put(6.8,2.3){\small $p$}
\pscircle*(8.5,1.2){0.07}
\pscircle*(8.5,2.4){0.07}
\pscircle*(10,2.52){0.07}
\psline(8.5,1.2)(11,1.2)
\psline(8.5,2.4)(11,2.4)
\psline(8.5,2.42)(11,2.8)
\psline(10,2.52)(11,2.52)
\put(11.2,1){\small $\lambda_2$}
\put(11.2,2){\small $\lambda_1$}
\put(11.2,2.45){\small $\lambda'_1$}
\put(11.2,3){\small $\lambda''_1$}
\end{picture}
\caption{Sliding $\lambda_1$ along $\lambda'_1$, and the auxiliary covering $p$}
\label{fig:arcslide}
\end{figure}
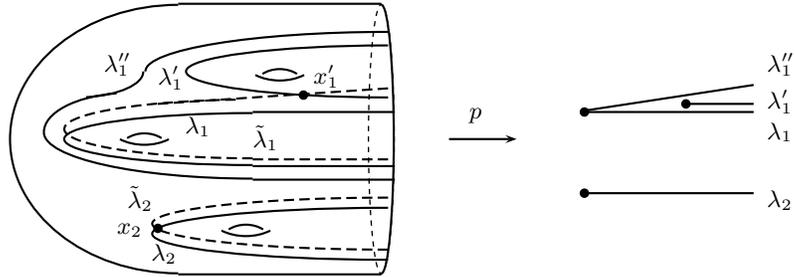

Consider $k+1$ disjoint properly embedded arcs
$\lambda_1,\dots,\lambda_k,\lambda'_1$ in $F$, with boundary in $\partial
F\setminus Z$, and such that an end point of $\lambda'_1$ lies immediately
after an end point of $\lambda_1$ along a component of $\partial F\setminus
Z$. Let $\lambda''_1$ be the arc obtained by sliding $\lambda_1$ along
$\lambda'_1$. Finally, denote by $\tilde\lambda_1,\dots,\tilde\lambda_k$
a collection of arcs obtained by slightly perturbing $\lambda_1,\dots,
\lambda_k$ in the positive direction, with each $\tilde\lambda_i$ intersecting 
$\lambda_i$ in a single point $x_i\in U$, and $\tilde\lambda_1$ intersecting 
$\lambda'_1$ in a single point $x'_1$ which lies near the boundary;
see Figure \ref{fig:arcslide}. Let $L=\lambda_1\times\dots\times \lambda_k$,
$L'=\lambda'_1\times\lambda_2\times\dots\times \lambda_k$, and
$L''=\lambda''_1\times\lambda_2\times\dots\times \lambda_k$.
Then the point $(x'_1,x_2,\dots,x_k)\in
(\tilde\lambda_1\times\dots\times \tilde\lambda_k)\cap (\lambda'_1\times\lambda_2
\times\dots\times \lambda_k)$ determines (via the appropriate continuation
map between Floer complexes, to account for the need to further perturb
$L$) an element of $\hom(L,L')$, which we call $u$. The following result
is essentially Lemma 5.2 of \cite{fuksymg}.

\begin{lemma}[\cite{fuksymg}]\label{l:arcslide}
In the $A_\infty$-category of twisted complexes $Tw\,\F(\Sym^k(F),Z)$, $L''$ is quasi-isomorphic to the mapping cone of
$u$.
\end{lemma}

The main idea is to consider an auxiliary simple branched covering
$p:\hat{F}\to\C$ for which the arcs $\lambda_1,\lambda'_1,\dots,\lambda_k$
are Lefschetz thimbles (i.e., lifts of half-lines), with the critical
value for $\lambda'_1$ lying immediately next to that for $\lambda_1$
and so that the monodromies at the corresponding critical values are
transpositions with one common index (see Figure \ref{fig:arcslide} right).
The objects $L,L',L''$ can be viewed as Lefschetz thimbles for 
the Lefschetz fibration induced by $p$ on the symmetric product;
in the corresponding Fukaya category, the statement that $L''\simeq
\mathrm{Cone}(u)$ follows from a general result of Seidel \cite[Proposition 18.23]{SeBook}.
The lemma then follows from exactness of the relevant acceleration functor.
See \S 5 of \cite{fuksymg} for details.

The other useful fact is that sliding one factor of $L$ over
another factor of $L$ only affects $L$ by a Hamiltonian isotopy.
For instance, in the above situation, $\lambda_1\times\lambda'_1\times
\lambda_3\times\dots\times \lambda_k$ and
$\lambda''_1\times\lambda'_1\times
\lambda_3\times\dots\times \lambda_k$ are Hamiltonian isotopic. (This is an
easy consequence of the main result in \cite{PerHH}).
\medskip

Returning to the collection of arcs $\underline\alpha'$ on the surface $F$, 
assume that $\alpha'_m$ can be erased without losing the property that every 
component of the complement is a disc carrying at most one point of $Z$. 
Then one of the connected components of $F\setminus (\alpha'_1\cup\dots\cup
\alpha'_m)$ is a disc $\Delta$ which contains no
point of $Z$, and whose boundary consists of portions of $\partial F$ and
the
arcs $\alpha'_m,\alpha'_{i_1},\dots,\alpha'_{i_r}$ (with $i_1,\dots,i_r$ distinct
from $m$, but not necessarily pairwise distinct) in that order. Then
the arc obtained by sliding $\alpha'_{i_1}$ successively over
$\alpha'_{i_2},\dots,\alpha'_{i_r}$ is isotopic to $\alpha'_m$. Hence,
by Lemma \ref{l:arcslide}, for $m\in s$ the object $D'_s$ 
can be expressed as a twisted complex built from the objects
$D'_{s_j}$, where $s_j=(s\cup\{i_j\})\setminus \{m\}$, for
$j\in \{1,\dots,r\}$ such that $i_j\not\in s$.

\section{Yoneda embedding and invariants of bordered 3-manifolds}
\label{s:yoneda}

Let $\FF=(F,Z,\underline\alpha)$ be a decorated surface, and assume
that every component of $F\setminus (\bigcup\alpha_i)$ is a disc carrying at
most one point of $Z$. By Theorem \ref{thm:generate} the partially
wrapped Fukaya category $\F(\Sym^k(F),Z)$ is generated by the product
objects $D_s$, $s\in \SS^n_k$. In fact Theorem \ref{thm:generate} continues
to hold if we consider the extended category $\F^\#(\Sym^k(F),Z)$ instead
of $\F(\Sym^k(F),Z)$; see
Proposition 6.3 of \cite{fuksymg}. (The key point is that the only
generalized Lagrangians we consider are compactly supported in the interior
of $\Sym^k(F)$, and Seidel's argument for generation of compact objects
by Lefschetz thimbles still applies to them.)

To each object $\T$ of $\F^\#(\Sym^k(F),Z)$ we can associate a right
$A_\infty$-module over the algebra $\A=\A(\FF,k)$,
$$\Y(\T)=\Y^r(\T)=\textstyle\bigoplus\limits_{s\in \SS^n_k}
\hom(\T,D_s)\in \text{mod-}\A,$$ where the module maps $m_\ell:
\Y(\T)\otimes \A^{\otimes (\ell-1)}\to \Y(\T)$ are defined by
products and higher compositions in $\F^\#(\Sym^k(F),Z)$.
Moreover, given two objects $\T_0,\T_1$, compositions in the partially 
wrapped Fukaya category yield a natural map from $\hom(\T_0,\T_1)$ 
to $\hom_{\text{mod-}\A}(\Y(\T_1),\Y(\T_0))$, as well as higher order
maps. Thus, we obtain a contravariant $A_\infty$-functor
$\Y:\F^\#(\Sym^k(F),Z)\to\text{mod-}\A(\FF,k)$: the
{\it right Yoneda embedding}.

\begin{proposition}\label{prop:embed}
Under the assumptions of Theorem \ref{thm:generate}, $\mathcal{Y}$ is a cohomologically full and faithful (contravariant)
embedding.
\end{proposition}

Indeed, the general Yoneda embedding into $\text{mod-}\F^\#(\Sym^k(F),Z)$
is cohomologically full and faithful (see e.g.\ \cite[Corollary
2.13]{SeBook}), while Theorem \ref{thm:generate} (or rather its analogue
for the extended Fukaya category) implies that the natural functor from
$\text{mod-}\F^\#(\Sym^k(F),Z))$ to $\text{mod-}\A(\FF,k)$ given by
restricting an arbitrary $A_\infty$-module to the subset of objects $\{D_s,\
s\in\SS^n_k\}$ is an equivalence.

We can similarly consider the {\it left Yoneda embedding}
to left $A_\infty$-modules over $\A(\FF,k)$, namely the (covariant) $A_\infty$-functor
$\Y^\ell:\F^\#(\Sym^k(F),Z)\to\A(\FF,k)\text{-mod}$ which sends the
object $\T$ to $\Y^\ell(\T)=\bigoplus_{s\in \SS^n_k} \hom(D_s,\T)$.

\begin{lemma}\label{l:Aop}
Denote by $-\FF=(-F,Z,\underline\alpha)$ the decorated surface obtained 
by orientation reversal. Then $\A(-\FF,k)$ is isomorphic to
the opposite $A_\infty$-algebra $\A(\FF,k)^{op}$.
\end{lemma}

\proof Given $s_0,\dots,s_\ell\in \SS^n_k$, and any positive perturbation
$(\tilde{D}_{s_0},\dots,\tilde{D}_{s_\ell})$ of the sequence $(D_{s_0},
\dots,D_{s_\ell})$ in $\Sym^k(F)$ relatively to $Z$, the reversed sequence
$(\tilde{D}_{s_\ell},\dots,\tilde{D}_{s_0})$ is a positive perturbation of
$(D_{s_\ell},\dots,D_{s_0})$ in $\Sym^k(-F)$.  Thus, the holomorphic
discs in $\Sym^k(-F)$ which contribute to the product operation
$m_\ell:\hom(D_{s_\ell},D_{s_{\ell-1}})\otimes
\dots\otimes\hom(D_{s_{1}},D_{s_0})\to \hom(D_{s_\ell},D_{s_0})$ in
$\A(-\FF,k)$ are exactly the complex conjugates of the holomorphic discs in
$\Sym^k(F)$ which contribute to $m_\ell:\hom(D_{s_0},D_{s_{1}})\otimes
\dots\otimes\hom(D_{s_{\ell-1}},D_{s_\ell})\to \hom(D_{s_0},D_{s_\ell})$ in
$\A(\FF,k)$.
\endproof

Hence, left $A_\infty$-modules over $\A=\A(\FF,k)$ can be interchangeably
viewed as right $A_\infty$-modules over $\A^{op}=\A(-\FF,k)$; more
specifically, given a generalized Lagrangian $\T$ in $\Sym^k(F)$ and
its conjugate $-\T$ in $\Sym^k(-F)$, the left Yoneda module
$\Y^\ell(\T)\in \A\text{-mod}$ is the same as the right Yoneda module
$\Y^r(-\T)\in \text{mod-}\A^{op}$.

Moreover, the left and right Yoneda embeddings are dual to each other:

\begin{lemma}\label{l:dual}
For any object $\T$, the modules $\Y^\ell(\T)\in\A\text{-mod}$
and $\Y^r(\T)\in\text{mod-}\A$ satisfy $\Y^r(\T)\simeq \hom_{\A\text{-mod}}
(\Y^\ell(\T),\A)$ and $\Y^\ell(\T)\simeq
\hom_{\text{mod-}\A}(\Y^r(\T),\A)$ (where $\A$ is viewed as an
$A_\infty$-bimodule over itself).
\end{lemma}

\proof By definition, $\Y^r(\T)=\hom(\T,\bigoplus_s D_s)$ (working in an
additive enlargement of $\F^\#(\Sym^k(F),Z)$), with the right $A_\infty$-module
structure coming from right composition (and higher products) with
endomorphisms of $\bigoplus_s D_s$. However, the left Yoneda
embedding functor is full and faithful,
and maps $\T$ to $\Y^\ell(\T)$ and $\bigoplus_s D_s$ to $\A$. Hence,
as chain complexes $\hom(\T,\bigoplus_s D_s)\simeq
\hom_{\A\text{-mod}}(\Y^\ell(\T),\A)$. Moreover this quasi-isomorphism
is compatible with the right module structures (by
functoriality of the left Yoneda embedding). The other statement is proved
similarly, by applying the right Yoneda functor (contravariant, full and
faithful) to prove that $\Y^\ell(\T)=\hom(\bigoplus_s D_s,\T)\simeq \hom_{\text{mod-}\A}
(\Y^r(\T),\A)$.\endproof

All the ingredients are now in place for the proof of
Theorem \ref{cor:pairing2} (and other similar pairing results). 
Consider as in the introduction a closed 3-manifold $Y$ obtained by 
gluing two 3-manifolds $Y_1,Y_2$ with $\partial Y_1=-\partial Y_2=
F\cup D^2$ along their common boundary, and equip the surface
$F$ with boundary
marked points $Z$ and a collection $\underline\alpha$ of disjoint 
properly embedded arcs such that
the decorated surface $\FF=(F,Z,\underline\alpha)$ satisfies the assumption
of Theorem \ref{thm:generate}. 

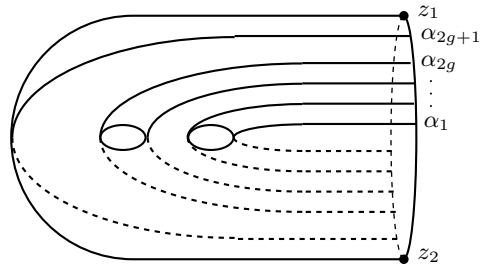
\begin{figure}[b]
\hfil
\setlength{\unitlength}{0.9cm}
\begin{picture}(5.5,3.5)(-0.5,-1.65)
\psset{unit=\unitlength}
\psellipticarc[linewidth=0.5pt,linestyle=dashed,dash=2pt
2pt](4.5,0)(0.2,1.8){90
}{-90}
\psellipticarc(4.5,0)(0.2,1.8){-90}{90}
\psline[linearc=1.8](4.5,1.8)(-1.3,1.8)(-1.3,-1.8)(4.5,-1.8)
\pscircle*(4.5,1.8){0.07} \put(4.7,1.8){\small $z_1$}
\pscircle*(4.5,-1.8){0.07} \put(4.7,-1.8){\small $z_2$}
\psellipse(0.35,0)(0.35,0.2)
\psellipse(1.65,0)(0.35,0.2)
\psellipticarc(3,0)(1.03,0.21){90}{180}
\psellipticarc(3,0)(1.71,0.51){90}{180}
\psellipticarc(3,0)(2.3,0.81){90}{180}
\psellipticarc(3,0)(3,1.11){90}{180}
\psellipticarc(2,0)(3.3,1.5){90}{180}
\psline(3,0.2)(4.68,0.2)
\psline(3,0.5)(4.65,0.5)
\psline(3,0.8)(4.65,0.8)
\psline(3,1.1)(4.6,1.1)
\psline(2,1.5)(4.6,1.5)
\psellipticarc[linestyle=dashed,dash=2pt 2pt](3,0)(1.03,0.215){180}{270}
\psellipticarc[linestyle=dashed,dash=2pt 2pt](3,0)(1.71,0.515){180}{270}
\psellipticarc[linestyle=dashed,dash=2pt 2pt](3,0)(2.3,0.815){180}{270}
\psellipticarc[linestyle=dashed,dash=2pt 2pt](3,0)(3,1.115){180}{270}
\psellipticarc[linestyle=dashed,dash=2pt 2pt](2,0)(3.3,1.5){180}{270}
\psline[linestyle=dashed,dash=2pt 2pt](3.05,-0.2)(4.32,-0.2)
\psline[linestyle=dashed,dash=2pt 2pt](3.05,-0.5)(4.32,-0.5)
\psline[linestyle=dashed,dash=2pt 2pt](3.05,-0.8)(4.35,-0.8)
\psline[linestyle=dashed,dash=2pt 2pt](3.05,-1.1)(4.38,-1.1)
\psline[linestyle=dashed,dash=2pt 2pt](2.05,-1.5)(4.4,-1.5)
\put(4.79,0.15){\small $\alpha_1$}
\put(4.75,1.05){\small $\alpha_{2g}$}
\put(4.75,1.45){\small $\alpha_{2g+1}$}
\psline[linestyle=dotted](4.9,0.85)(4.9,0.45)
\end{picture}
\caption{Decorating $F$ with two marked points and $2g+1$ arcs}
\label{fig:2g+1}
\end{figure}

\begin{remark*}
The natural choice in view
of Lipshitz-Ozsv\'ath-Thurston's work on bordered Heegaard-Floer
homology \cite{LOT} is to equip $F$ with 
a single marked point and a collection
of $2g$ arcs that decompose it into a single disc, e.g.\ as in Figure
\ref{fig:2g}. However, 
one could also equip $F$ with {\it two} boundary marked points and
$2g+1$ arcs, by viewing $F$ as a
double cover of the unit disc with $2g+1$ branch points and proceeding
as in \S \ref{ss:332}; see Figure \ref{fig:2g+1}. While this yields a
larger generating set, with $\binom{2g+1}{k}$ objects instead of 
$\binom{2g}{k}$, the resulting algebra remains combinatorial in nature
(by Proposition \ref{prop:nice}) and it is more familiar from the
perspective of symplectic geometry, since we
are now dealing with the Fukaya category of a Lefschetz fibration on
the symmetric product. Among other nice features, the generators are exceptional
objects, and the algebra is directed.
\end{remark*}

Let us now return to our main argument. As explained in \S \ref{ss:LP},
the work of Lekili and Perutz \cite{LP} associates to the 3-manifolds
$Y_1$ and $-Y_2$ (viewed as sutured cobordisms from $D^2$ to $F$)
two generalized Lagrangian submanifolds 
$\T_{Y_1}$ and $\T_{-Y_2}$ of $\Sym^g(F)$, with the property
that $\widehat{CF}(Y)$ is quasi-isomorphic to $\hom_{\F^\#(\Sym^g(F))}(
\T_{Y_1},\T_{-Y_2})$. However, by Proposition \ref{prop:embed} we have
$$\hom_{\F^\#(\Sym^g(F))}(\T_{Y_1},\T_{-Y_2})\simeq
\hom_{\text{mod-}\A(\FF,g)}(\Y(\T_{-Y_2}),\Y(\T_{Y_1}))$$
where $\Y=\Y^r$ denotes the right Yoneda embedding functor.
Moreover, using Lemma \ref{l:dual}, and setting $\A=\A(\FF,g)$, we have:
\begin{eqnarray*}
\Y^r(\T_{Y_1})\otimes_{\A}\Y^\ell(\T_{-Y_2})&\simeq&
\Y^r(\T_{Y_1})\otimes_{\A}\hom_{\text{mod-}\A}(\Y^r(\T_{-Y_2}),\A)\\
&\simeq& \hom_{\text{mod-}\A}(\Y^r(\T_{-Y_2}),\Y^r(\T_{Y_1})\otimes_{\A}\A)\\
&\simeq& \hom_{\text{mod-}\A}(\Y^r(\T_{-Y_2}),\Y^r(\T_{Y_1})).
\end{eqnarray*}
Finally, by the discussion after Lemma \ref{l:Aop}, we can identify
the left $\A(\FF,g)$-module $\Y^\ell(\T_{-Y_2})$ with the right module
$\Y^r(\T_{Y_2})\in \text{mod-}\A(-\FF,g)$. This completes the proof of
Theorem \ref{cor:pairing2}.
\medskip

Turning to the case of more general cobordisms, recall
that the construction of Lekili and Perutz associates 
to a sutured manifold $Y$ with $\partial Y=(-F_-)\cup F_+$ a
generalized Lagrangian correspondence $\T_Y$ from $\Sym^{k_-}(F_-)$
to $\Sym^{k_+}(F_+)$ (where $k_+-k_-=g(F_+)-g(F_-)$), i.e.\ an object of
$\F^\#(\Sym^{k_-}(-F_-)\times\Sym^{k_+}(F_+))$.

Equip the surfaces $F_-$ and $F_+$ with sets of boundary marked points 
$Z_\pm$ and two collections 
$\underline{\alpha}_\pm$ of properly
embedded arcs such that the decorated surfaces $\FF_\pm=
(F_\pm,Z_\pm,\underline\alpha_\pm)$ satisfy the assumption of Theorem~\ref{thm:generate}.
Considering products of $k_\pm$ of the arcs in $\underline\alpha_\pm$, we
have two collections of product Lagrangian submanifolds $D_{\pm,s}$
($s\in\SS_\pm$) in $\Sym^{k_\pm}(F_\pm)$. By a straightforward
generalization of 
Theorem \ref{thm:generate}, the partially wrapped
Fukaya category $\F^\#(\Sym^{k_-}(-F_-)\times\Sym^{k_+}(F_+),Z_-\sqcup Z_+)$ 
is generated by the product objects $(-D_{-,s})\times D_{+,t}$ for $(s,t)\in
\SS_-\times\SS_+$. 
Indeed,
$\Sym^{k_-}(-F_-)\times\Sym^{k_+}(F_+)$ is a connected component of 
$\Sym^{k_-+k_+}((-F_-)\sqcup F_+)$,
and the proof of Theorem~\ref{thm:generate} applies without modification
to the 
disconnected decorated surface $(-\FF_-)\sqcup \FF_+=((-F_-)\sqcup F_+,
Z_-\sqcup Z_+,\underline\alpha_-\sqcup \underline\alpha_+)$.
Hence, as before, the Yoneda construction
$$\textstyle \Y(\T_Y)=\bigoplus\limits_{(s,t)\in \SS_-\times\SS_+}
\hom(\T_Y,(-D_{-,s})\times D_{+,t})$$ defines a cohomologically full
and faithful embedding into the category of right $A_\infty$-bimodules
over $\A(-\FF_-,k_-)$ and $\A(\FF_+,k_+)$, or equivalently, the
category of $A_\infty$-bimodules $\A(\FF_-,k_-)\text{-mod-}\A(\FF_+,k_+)$.
This property is the key ingredient that makes it possible to relate 
compositions of generalized Lagrangian correspondences (i.e., gluing of sutured
cobordisms) to algebraic operations on $A_\infty$-bimodules, as in 
Conjecture \ref{conj:pairing} for instance.

\section{Relation to bordered Heegaard-Floer homology}
\label{s:bordered}

Consider a sutured 3-manifold $Y$, with 
$\partial Y=(-F_-)\cup (\Gamma\times[0,1])\cup
F_+$, and pick decorations $\FF_\pm=(F_\pm,Z_\pm,\underline\alpha_\pm)$ of
$F_\pm$. Assume for simplicity that $Z_+=Z_-$.
Denote by $g_\pm$ the genus of $F_\pm$, and
by $n_\pm$ the number of arcs in $\underline\alpha_\pm$. Choose a Morse
function $f:Y\to [0,1]$ with index 1 and 2 critical points only,
such that $f^{-1}(1)=F_-$ and $f^{-1}(0)=F_+$. Assume that all the 
index 1 critical points lie in $f^{-1}((0,\frac12))$ and all the
index 2 critical points lie in $f^{-1}((\frac12,1))$. Also pick
a gradient-like vector field for $f$, tangent to the boundary along
$\Gamma\times[0,1]$, and equip the level sets of $f$ with
complex structures such that the gradient flow induces biholomorphisms away
from the critical locus.
The above data determines a {\em bordered Heegaard diagram}
on the surface $\Sigma=f^{-1}(\frac12)$ of genus $\bar{g}=g(\Sigma)$, consisting of:
\begin{itemize}
\item
$\bar{g}-g_+$ simple closed curves
$\alpha_1^c,\dots,\alpha_{\bar{g}-g_+}^c$, where $\alpha_i^c$ is the set of
points of
$\Sigma$ from which the downwards gradient flow converges to the $i$-th index 1 critical
point;
\item
$n_+$ properly embedded arcs $\alpha_1^a,\dots,\alpha_{n_+}^a$, where
$\alpha_i^a$ is the set of points of $\Sigma$
from which the downwards gradient flow ends at a point of $\alpha_{+,i}\subset
F_+$;
\item
$\bar{g}-g_-$ simple closed curves
$\beta_1^c,\dots,\beta_{\bar{g}-g_-}^c$, where $\beta_i^c$ is the set of
points of
$\Sigma$ from which the upwards gradient flow converges to the $i$-th index
2 critical point;
\item
$n_-$ properly embedded arcs $\beta_1^a,\dots,\beta_{n_-}^a$, where
$\beta_i^a$ is the set of points of $\Sigma$
from which the upwards gradient flow ends at a point of $\alpha_{-,i}\subset
F_-$;
\item a finite set $Z$ of boundary marked points 
(which match with $Z_\pm$ under the gradient flow).
\end{itemize}

Given integers $\bar{k},k_+,k_-$ satisfying
$\bar{k}-\bar{g}=k_+-g_+=k_--g_-$, we can view the generalized Lagrangian
correspondence $\T_Y$ associated to $Y$ as the composition of the correspondence
$T_\beta\subset\Sym^{k_-}(-F_-)\times\Sym^{\bar{k}}(\Sigma)$ determined by
$f^{-1}([\frac12,1])$ and the correspondence $T_\alpha\subset
\Sym^{\bar{k}}(-\Sigma)\times \Sym^{k_+}(F_+)$ determined by $f^{-1}([0,
\frac12])$.\smallskip

The $A_\infty$-bimodule $\Y(\T_Y)\in \A(\FF_-,k_-)\text{-mod-}\A(\FF_+,k_+)$
associated to $Y$ can then be understood entirely in terms of the symmetric
product $\Sym^{\bar{k}}(\Sigma)$. Namely, denote by 
$\bar\F^\#=\bar\F^\#(\Sym^{\bar{k}}(\Sigma),Z)$ a partially wrapped Fukaya category
defined similarly to the construction in \S \ref{s:fukaya}, except we also allow 
objects which are products of mutually disjoint simple closed
curves and properly embedded arcs in $\Sigma$. 

The Lagrangian correspondences $-T_\alpha$ and $T_\beta$
induce $A_\infty$-functors $\Phi_{\alpha}$ and~$\Phi_\beta$ from 
$\F^\#(\Sym^{k_\pm}(F_\pm),Z_\pm)$ to $\bar\F^\#$.
Considering the product Lagrangians $D_{-,s}$ for
\hbox{$s\in \SS_-=\SS^{n_-}_{k_-}$}, the description of the geometry of the 
correspondence $T_\beta$ away from the 
diagonal~\cite{Perutz} (or the result of \cite{LP}) implies that
$\Phi_{\beta}(D_{-,s})$, i.e., the composition of $D_{-,s}$ with
the correspondence $T_\beta$, is Hamiltonian isotopic to 
$$\textstyle\Delta_{\beta,s}=\prod\limits_{i\in s}\beta_i^a\times
\prod\limits_{j=1}^{\ \bar{g}-g_-}\beta_j^c\subset \Sym^{\bar{k}}(\Sigma).$$
Similarly, for $t\in \SS_+=\SS^{n_+}_{k_+}$ the image of $D_{+,t}$ under the
correspondence $(-T_\alpha)$ is Hamiltonian isotopic to the product
$$\textstyle\Delta_{\alpha,t}=\prod\limits_{i\in t
}\alpha_i^a\times
\prod\limits_{j=1}^{\ \bar{g}-g_+}\alpha_j^c\subset \Sym^{\bar{k}}(\Sigma).$$
This implies the following result:

\begin{proposition}\label{prop:adjoint}
The $A_\infty$-bimodule $\Y(\T_Y)\in\A(\FF_-,k_-)\text{-mod-}\A(\FF_+,k_+)$
is quasi-isomorphic to $\bigoplus_{s,t}
\hom_{\bar\F^\#}(\Delta_{\beta,s},\Delta_{\alpha,t})$.
\end{proposition}

To clarify this statement, observe that $\Phi_\alpha$ induces
an $A_\infty$-homomorphism from $\A(\FF_+,k_+)=\bigoplus_{s,t}
\hom(D_{+,s},D_{+,t})$ to
$\A_\alpha=\bigoplus_{s,t} \hom_{\bar\F^\#}(\Delta_{\alpha,s},
\Delta_{\alpha,t})$. In fact, suitable choices in the construction
ensure that $\A_\alpha\simeq \A(\FF_+,k_+)\otimes H^*(T^{\bar{g}-g_+},\Z_2)$
and the map from $\A(\FF_+,k_+)$ to $\A_\alpha$ is simply given by
$x\mapsto x\otimes 1$. In any case, via $\Phi_\alpha$ we can view any right
$A_\infty$-module over $\A_\alpha$ as a right $A_\infty$-module over
$\A(\FF_+,k_+)$. 
Similarly, $\Phi_\beta$ induces an $A_\infty$-homomorphism
from $\A(\FF_-,k_-)$ to $\A_\beta=\bigoplus_{s,t}\hom_{\bar\F^\#}
(\Delta_{\beta,s},\Delta_{\beta,t})$, through which any left
$A_\infty$-module over $\A_\beta$ can be viewed as a left $A_\infty$-module
over $\A(\FF_-,k_-)$.

With this understood, Proposition \ref{prop:adjoint} essentially
follows from the fact that the $A_\infty$-functors induced by the
correspondences $T_\alpha$ and $(-T_\alpha)$ on one hand, and
$T_\beta$ and $(-T_\beta)$ on the other hand, are adjoint to each other;
see Proposition 6.6 in \cite{fuksymg} for the case of $A_\infty$-modules.
\medskip

The case where one of $k_\pm$ vanishes, say $k_-=0$, is of particular
interest; then the $\beta$-arcs play no role whatsoever, and we only need to consider
the product torus $T_\beta=\beta_1^c\times\dots\times \beta^c_{\bar{g}-g_-}
\subset \Sym^{\bar{k}}(\Sigma)$. This happens for instance when
$F_-$ is a disc, i.e.\ when $Y$ is a 3-manifold with boundary
$\partial Y=F_+\cup D^2$ viewed as a sutured cobordism from $D^2$ to $F_+$. 
(This corresponds to the situation considered in \cite{LOT}; in this case
we have $\bar{k}=\bar{g}$ and $k_+=g_+$). 

In this situation, the statement of Proposition
\ref{prop:adjoint} becomes that the right $A_\infty$-module
$\Y(\T_Y)\in\text{mod-}\A(\FF_+,k_+)$ is quasi-isomorphic to
$\bigoplus_{t\in \SS_+\!} \hom_{\bar\F^\#}(T_\beta,\Delta_{\alpha,t})$. 
Then we have the following result
(Proposition 6.5 of~\cite{fuksymg}):

\begin{proposition}\label{prop:CFA}
The right $A_\infty$-modules over $\A(\FF_+,k_+)$ constructed by
Yoneda embedding, $\Y(\T_Y)\simeq \bigoplus_{t\in \SS_+}
\hom_{\bar\F^\#}(T_\beta,\Delta_{\alpha,t})$, and by bordered
Heegaard-Floer homology, $\smash{\widehat{CFA}}(Y)$, are quasi-isomorphic.
\end{proposition}

The fact that $\bigoplus_{t}\hom_{\bar\F^\#}(T_\beta,\Delta_{\alpha,t})$
and $\widehat{CFA}(Y)$ are quasi-isomorphic (in fact isomorphic) as chain complexes is a
straightforward consequence of the definitions. Comparing the
module structures requires a comparison of the moduli
spaces of holomorphic curves which determine the
module maps; this can be done via a neck-stretching argument, see \cite[Proposition
6.5]{fuksymg}.

\begin{remark*}
Another special case worth mentioning is when $k_+=k_-=0$, which requires
the sutured manifold $Y$ to be balanced in the sense of \cite{juhasz}.
Then we can discard all the arcs from the Heegaard diagram, and
$\Y(\T_Y)\simeq \hom_{\bar\F^\#}(T_\beta,T_\alpha)$ is simply 
the chain complex which defines the sutured Floer homology of \cite{juhasz}.
In this sense, bordered Heegaard-Floer homology and our constructions can 
be viewed as natural generalizations of Juh\'asz's sutured Floer homology.
(An even greater level of generality is considered in \cite{Zarev}.)
\end{remark*}

In light of the relation between $\Y(\T_Y)$ and $\widehat{CFA}(Y)$,
it is interesting to compare Theorem \ref{cor:pairing2} with the
pairing theorem obtained by Lipshitz, Ozsv\'ath and Thurston
for bordered Heegaard-Floer homology \cite{LOT}. In particular,
a side-by-side comparison suggests that the modules $\widehat{CFA}(Y)$ and
$\widehat{CFD}(Y)$ might be quasi-isomorphic.
\medskip

Another surprising aspect, about which we can only offer speculation,
is the seemingly different manners in which bimodules arise in the two
stories.  In our case, bimodules arise from sutured 3-manifolds viewed
as cobordisms between decorated surfaces, i.e.\ from bordered Heegaard
diagrams where both $\alpha$- and $\beta$-arcs are simultaneously
present; and pairing results arise from ``top-to-bottom'' stacking of
cobordisms. On the other hand, the work of Lipshitz, Ozsv\'ath and 
Thurston~\hbox{\cite{LOT,LOT2}} provides a different construction of bimodules associated
to cobordisms between decorated surfaces, involving diagrams in which there
are no $\beta$-arcs; and pairing results arise from ``side-by-side''
gluing of bordered Heegaard diagrams.

As a possible way to understand ``side-by-side'' gluing in our framework,
observe that given two decorated surfaces
$\FF_i=(F_i,Z_i,\underline\alpha_i)$ for $i=1,2$, and given two
points $z_1\in Z_1$ and $z_2\in Z_2$, we can form the {\it boundary
connected sum} $F=F_1\cup_\partial F_2$ of $F_1$ and $F_2$ by attaching a 1-handle (i.e., a
band) to small
intervals of $\partial F_1$ and $\partial F_2$ containing $z_1$ and $z_2$
respectively. The surface $F$ can be equipped with the set of
marked points $Z=(Z_1\setminus \{z_1\})\cup (Z_2\setminus \{z_2\})\cup
\{z_-,z_+\}$, where $z_-$ and $z_+$ lie on either side of the
connecting handle, and the collection of properly embedded arcs
$\underline\alpha=\underline\alpha_1\cup\underline\alpha_2$.
Assume moreover that $\FF_1$ and $\FF_2$ satisfy the conditions of 
Proposition \ref{prop:nice}, so that the associated algebras are
honest differential algebras. 
Denoting by $\FF$ the decorated surface $(F,Z,\underline\alpha)$,
it is then easy to check that $\A(\FF,k)\simeq \bigoplus\limits_{k_1+k_2=k}
\A(\FF_1,k_1)\otimes \A(\FF_2,k_2)$.

Now, given two 3-manifolds $Y_1,Y_2$ with boundary $\partial Y_i\simeq
F_i\cup D^2$, we can form their boundary connected sum
$Y=Y_1\cup_\partial Y_2$ by attaching a
1-handle at the points $z_1,z_2$; then $\partial Y=F\cup D^2$, and the
bordered Heegaard diagram representing $Y$ is
simply the boundary connected sum of the 
bordered Heegaard diagrams representing $Y_1$ and $Y_2$.
Accordingly, the right $A_\infty$-module associated
to $Y$ is the tensor product (over the ground field $\Z_2$!) of the
right $A_\infty$-modules associated to $Y_1$ and~$Y_2$.
In the case where $\FF_1\simeq -\FF_2$, we can glue a standard handlebody
to $Y$ in order to obtain a closed 3-manifold $\bar{Y}$, namely the
result of gluing $Y_1$ and $Y_2$ along their entire boundaries
rather than just at small discs near the points $z_1,z_2$.
However, because the decorated surface $\FF$ never satisfies the
assumption of Theorem \ref{thm:generate} (the two new marked points
$z_\pm$ lie in the same component), the Yoneda functor to $A_\infty$-modules
over $\A(\FF,g)$ is not guaranteed to be full and faithful, so our
gluing result does not apply.

\end{document}